\documentclass[
3p,
times
]{elsarticle}

\usepackage{graphicx,color}

\usepackage[colorlinks=true, breaklinks=true, pdfstartview=FitV, linkcolor=red, citecolor=blue, urlcolor=black]
{hyperref}

\usepackage[T1]{fontenc}

\usepackage{amsmath}

\renewcommand{\Vec}[1]{\mbox{\boldmath$#1$}}

\journal{Commun Nonlinear Sci Numer Simulat} 

\begin{document}

\begin{frontmatter}

\title{Exploring the long-term dynamics of perturbed Keplerian motion in high degree potential fields$^{1}$}

\author[RA]{Martin Lara\corref{cor}\fnref{footnote1,footnote2}}
\ead{martin.lara@unirioja.es}

\author[CIBIR]{Rosario L\'opez\fnref{footnote1}}
\ead{rosario.lopez@unirioja.es}

\author[UR]{Iv\'an L.~P\'erez\fnref{footnote1}}
\ead{ivan.perez@unirioja.es}

\author[RA]{Juan F.~San-Juan\fnref{footnote1}}
\ead{juanfelix.sanjuan@unirioja.es}

\address[RA]{Dpto. Matematicas y Computacion, C/ Madre de Dios, 53, ES-26006 Logro\~no, Spain}
\address[CIBIR]{Center for Biomedical Research of La Rioja, C/ Piqueras 98, ES-26006 Logro\~no, Spain}
\address[UR]{Dep.~of Mechanical Engineering, University of La Rioja, C/ Luis de Ulloa, s.n., ES-26004 Logro\~no, Spain}
\fntext[footnote3]{Preliminary results were presented in the 7th International Conference on Astrodynamics Tool and Techniques, 6-9 November 2018, DLR Oberpfaffenhofen, Germany }

\cortext[cor]{Corresponding author }

\fntext[footnote1]{GRUCACI, University of La Rioja}
\fntext[footnote2]{Space Dynamics Group, Technical University of Madrid -- UPM}

\date{draft of \color{red}\today}

\begin{abstract}
The long-term dynamics of perturbed Keplerian motion is usually analyzed in simplified models as part of the preliminary design of artificial satellites missions. It is commonly approached by averaging procedures that deal with literal expressions in expanded form. However, there are cases in which the correct description of the dynamics may require full, contrary to simplified, potential models, as is, for instance, the case of low-altitude, high-inclination lunar orbits. In these cases, dealing with literal expressions is yet possible with the help of modern symbolic algebra systems, for which memory handling is no longer an issue. Still, the efficient evaluation of the averaged expressions related to a high fidelity potential is often jeopardized for the expanded character of the output of the automatic algebraic process, which unavoidably provides huge expressions that commonly comprise tens of thousands of literal terms. Rearrangement of the output to generate an efficient numerical code may solve the problem, but automatization of this kind of post-processing is a non trivial task due to the ad-hoc heuristic simplification procedures involved in the optimization process. However, in those cases in which the coupling of different perturbations is not of relevance for the analysis, the averaging procedure may preserve the main features of the structure of the potential model, thus avoiding the need of the typical blind computer-based brut force perturbation approach. Indeed, we show how standard recursions in the literature may be used to efficiently replace the brut force approach, in this way avoiding the need of further simplification to improve performance evaluation. In particular, Kaula's seminal recursion formulas for the gravity potential reveal clearly superior to the use of both expanded expressions and other recursions more recently proposed in the literature.
\end{abstract}


\end{frontmatter}

\section{Introduction}

The complexity of the orbital problems that aerospace engineers commonly confront makes that analytical solutions to them remain generally unknown. Because of that, simplified models that may capture the bulk of the dynamics of concern for a particular mission are customarily used in the preliminary steps of mission design \cite{Kozai1963,CuttingFrautnickBorn1978,ScheeresGumanVillac2001,LaraSanJuan2005,ColomboLuckingMcInnes2012}, but also in the search for efficient maneuvers for the end-of-life disposal of the spacecraft \cite{ArmellinSanJuanLara2015,Colomboetal2015}. Some of these simplified models have only two degrees of freedom (DOF) and, therefore, the phase space description can be approached with the usual tools of nonlinear dynamics, as the computation of Poincar\'e surfaces of section \cite{Broucke1994,JorbaMasdemont1999}, the numerical continuation of periodic orbits and other invariant manifolds \cite{Lara1999,GomezMondelo2001,Lara2003,LaraRussellVillac2007,RussellLara2007}, or the computation of different stability indicators which may disclose the existence of stable regions in real models, as well as the dynamical channels connecting distant ones \cite{LaraRussellVillac2006,LaraRussellVillac2007Meccanica,Villac2008}. Alternatively, an analytical process in which the dimension of the system is reduced to a 1-DOF integrable approximation of the original problem is customarily done by filtering the higher frequencies of the motion using perturbation theory. The integrable approximation is valid in some region of the phase space and reflects the main characteristics of the long-term dynamics, which can be explored without need of numerical integration by the (almost) instant representation of eccentricity-vector diagrams and inclination-eccentricity curves of frozen orbits \cite{CoffeyDepritDeprit1994,SanJuanLaraFerrer2006}. 
\par

There are, however, some cases in which approximate descriptions of dynamics, even qualitative ones, cannot be characterized by a simple model and, on the contrary, require to deal with full models. The paradigm of these problems is the design of low lunar orbits, a procedure in which, due to the irregular distribution of the mass of the moon, it is commonly accepted that, at least, a $50\times0$ truncation of the Selenepotential is required in the correct description of the long-term dynamics \cite{Konoplivetal1994,Roncoli2005,LaraSaedeleerFerrer2009,Lara2011}. This is not the unique case in which full models are required, and the use of higher order truncations of the gravity potential has also been  encouraged in the preliminary design of low-eccentricity, frozen earth orbits \cite{Cook1966,RosboroughOcampo1991,Cook1992}.
\par

The reduction to a 1-DOF system is still possible when dealing with full models, but the literal expressions to handle in the perturbation approach are enormous in this case, and, therefore, the assistance of a computer algebra system becomes essential \cite{Lara2010ICATT}. The current computational power makes that memory allocation is no longer a main concern, and literal expressions are customarily handled in expanded form up to higher degrees. However, rendering eccentricity-inclination diagrams and inclination-eccentricity curves of frozen orbits can no longer be considered an instant operation because the evaluation process is jeopardized by the unwieldy expressions obtained as the rough output of the symbolic processor. Indeed, series of tens of thousands, or even hundreds of thousands of literal terms are reported in the literature in relation to full perturbation models \cite{CoffeyNealSegermanTravisano1995}. A partial solution to this problem comes out from a post-processing of the crude symbolic output to turn it into a smartly rearranged sequence. However, in spite of broad simplification rules are part of the gunnery with which general symbolic algebra systems are equipped, designing efficient heuristic simplification procedures, which are required when dealing with implicit integration to avoid the limitations of series expansions in the eccentricity, is a quite challenging task because these kinds of simplifications strongly depend on the algebraic structure of particular problems. Because of that, mathematical simplification is still regarded as an art \cite{Knuth1997}, and the implementation of useful simplification rules definitely requires high skills in symbolic algebra manipulation jointly with an uncommon gift for identifying useful intermediate expressions. Hence, the development of efficient simplification and evaluation rules for coping with perturbed Keplerian motion remains only within the reach of a handful of selected experts \cite{CoffeyDeprit1980,Deprit1982,HealyTravisano1998}.
\par

Handling huge literal expressions cannot be avoided when higher orders of the perturbation approach are essential in the description of the dynamics \cite{WnukBreiter1991,Metrisetal1993,WnukJopek1994,LaraPalacianRussell2010,LaraPerezLopez2017}. On the contrary, in those cases in which the coupling of different perturbations is not of relevance in the analysis, the averaging procedure may preserve the main structural features of the potential model, thus avoiding the need of the typical blind computer-based, brut force approach. This is in particular the case of the design of low lunar orbits, whose dynamics is mainly driven by the Selenopotential, and provides the principal motivation for the current research. Indeed, analytical models based on low degree truncations of the Selenopotential not only fail in describing the quantitative behavior of high inclination, low lunar orbits, but they may provide completely wrong qualitative descriptions of the phase space. This fact is clearly illustrated in Fig.~(3) of \cite{LaraSaedeleerFerrer2009}, where it is shown that the actual argument of the periapsis of almost polar frozen orbits lies exactly in the opposite direction of the one predicted by a 7th degree truncation of the potential. Hence the convenience of preliminary exploration of the sensitivity of the dynamics, based on a gradual increase of the complexity of the gravitational potential in a harmonic coefficient-by-harmonic coefficient basis \cite{LaraFerrerSaedeleer2009,Lara2018Stardust}. Therefore, having available efficient procedures for evaluating algebraic expressions with a similar structure to the usual gravitational potential expansion is more than welcome.
\par

While the evaluation of the gravitational potential is efficiently done in Cartesian coordinates \cite{Cunningham1970,Pines1973,LundbergSchutz1988,FantinoCasotto2009}, the analytical investigations of the orbital elements evolution are better performed using Lagrange planetary equations, which require the reformulation of the gravitational potential in orbital elements. The formulation of the gravity potential  in orbital elements shows the contribution of secular, long-, and short-period corrections to the pure Keplerian motion \cite{Kozai1959}. Relevant aspects of the dynamics are disclosed when focusing on the effects of secular and long-period terms, which can be studied after averaging the contribution of short-period effects, namely, those related to the period in which the mean anomaly advances by $2\pi$. From the mathematical point of view, this averaging is the result of a transformation from osculating to ``mean'' elements, whose explicit computation is of the utmost importance when dealing with unstable dynamics, as is the common case of mapping orbits about planetary satellites \cite{Lara2008,LaraPalacianYanguasCorral2010}.
\par

The direct conversion of the gravitational potential into orbital elements is quite feasible for low degree and order truncations. On the other hand, the computational burden grows enormously when converting higher orders, yet the use of computer algebra systems makes the task possible. Alternatively, efficient recursion formulas for the conversion of the gravitational potential into orbital elements up to arbitrary degree and order have been provided by Kaula \cite{Kaula1961,Kaula1966}, who organizes the potential as a multivariate Fourier series. In these series, the arguments of the trigonometric functions involve linear combinations of the mean anomaly, the argument of the periapsis, and the right ascension of the ascending node, as well as the Greenwich angle, whereas the coefficients of the trigonometric terms are obtained from Kaula's inclination and eccentricity functions ---the latter being special cases of Hansen coefficients \cite{Hansen1855}--- as well as the different powers of the inverse of the orbit semi-major axis. Equinoctial elements \cite{ArsenaultFordKoskela1970} are commonly used to avoid the singularities of circular and equatorial orbits that are inherent to the classical Keplerian elements formulation. Formulas for the conversion of the gravitational potential into equinoctial elements are provided in \cite{BrouckeCefola1972,CefolaBroucke1975}. 
\par

Existing recursions are not limited to the construction of the gravitational potential alone, and are also available for the secular and long-period terms of the gravitational potential that define the mean elements dynamics. Indeed, when tesseral resonances are not of concern, the mean elements potential is reduced to the zonal harmonics contribution, a case in which Kaula's eccentricity functions are exact, contrary to the expansions in the eccentricity required in other instances. For this case, equivalent recursions to those of Kaula, but taking the point of view of perturbation theory using Delaunay canonical variables, have been later made available and complemented with analogous recursions for the construction of the generating function of the canonical transformation that performs the short-period averaging of the zonal potential, yet limited to a first order perturbation approach \cite{Saedeleer2005}. Analogous, but more involved recursions, exist also for third-body perturbations in different simplified models \cite{Kaula1962,LaskarBoue2010,PalacianVanegasYanguas2017}.
\par

The use of recursion notably speed evaluation of full gravitational fields, and slight modifications of the recursions in \cite{Saedeleer2005} together with additional recursions for the frozen orbits equation ease rendering phase space diagrams in real time \cite{LaraSaedeleerFerrer2009} even in the case in which a 100th degree truncation of the zonal potential is taken into account \cite{Lara2011}. Even though, we will show that the performance of the recursions in \cite{Saedeleer2005} is clearly exceeded by Kaula's original recursion formulas for the averaged potential. Note, however, that these kinds of recursions are constrained to first order effects in the perturbation approach. That is, they are useful only in those cases in which there is not a clear prevalence of any of the potential terms, which, therefore, are all taken to be of the same order in the perturbation arrangement. The moon and Venus are relevant instances of this kind of gravitational potential. On the contrary, it is well known that the dynamics about earth-like bodies is dominated by the zonal harmonic coefficient of the second degree ($C_{2,0}$), whereas the rest of coefficients are $\mathcal{O}(C_{2,0}^2)$. In consequence, second order effects of $C_{2,0}$ are as much important as first order effects of the other coefficients, and, because of that, they must be taken into account in studies of the long-term dynamics about earth-like bodies \cite{Brouwer1959,Garfinkel1959,Kozai1959,CoffeyDepritDeprit1994}. Hence, to broaden applicability of Kaula's recursions to this common case, they are complemented with the terms that model the long-term second order effects due to the zonal harmonic of the second degree. As far as recursions for a second order averaging are not yet available, we simply add the expanded expressions of this second order effect to the (first order) recursion formulas. 
\par


The paper is organized as follows. We present basic facts on the formulation of the gravitational potential in orbital elements in Section \ref{s:zonal}, were we also recall fundamental formulas that are pertinent to approach the zonal potential Hamiltonian by perturbations taking benefit of Kaula's recursions. To properly deal with second order effects of the zonal harmonic of the second degree in those cases in which they are relevant, the elimination of the parallax \cite{Deprit1981,LaraSanJuanLopezOchoa2013b,LaraSanJuanLopezOchoa2013c} is carried out first, in Section \ref{s:parallax}, where it is shown that this simplification does not affect the basic structure of Kaula's recursions. Kaula's eccentricity functions are recovered at this stage, which can be expressed in closed form because we only deal with the zonal part of the gravitational potential. Remaining short-period effects after the elimination of the parallax are removed in Section \ref{s:delaunay} by performing a Delaunay normalization \cite{Deprit1982}, leading to the compact formulation of a mean elements Hamiltonian by means Kaula-type recursions, which are supplemented with the expanded terms corresponding to the second order effects of the zonal harmonic coefficient of the second degree. Performance comparisons between Kaula and other recursions in the literature and presented in Section \ref{s:performance}, showing the superiority of Kaula's formulation for higher degree truncations of the mean elements zonal potential. Finally, the feasibility of displaying eccentricity-vector diagrams and inclination-eccentricity curves in real time, without limiting to the case of low eccentricities or linearized equations \cite{Cook1966,RosboroughOcampo1991,Cook1992}, is illustrated in Section \ref{s:averagedflow}, where a coefficient-by-coefficient approach is used to ascertain the correct truncation of the Selenopotential for increasing altitudes of a lunar orbiter over the surface of the moon. The different transformation and simplifications carried out to obtain the long-term dynamics are based on Deprit's implementation of the Lie transforms method \cite{Deprit1969}. For completeness, the basics of this perturbation method, which is standard these days \cite{BoccalettiPucacco1998v2,MeyerHall1992}, are summarized in \ref{s:LieTransforms}.
\par

\section{The Hamiltonian of the zonal potential} \label{s:zonal}

Solution of the Laplace's equation in spherical coordinates leads to the usual expansion of the gravity potential \cite[see][for instance]{MacMillan1958}
\begin{equation} \label{geopot}
U=-\frac{\mu}{r}\sum_{n\ge{0}}\frac{R_\oplus^n}{r^n}
\sum_{m=0}^n\left(C_{n,m}\cos{m\lambda}+S_{n,m}\sin{m\lambda}\right)P_{n,m}(\sin\varphi),
\end{equation}
where $r$, $\varphi$, and $\lambda$ stand for radius, geocentric latitude and longitude, respectively;  $P_{n,m}$ are associated Legendre polynomials, and $C_{n,m}$ and $S_{n,m}$ are harmonic coefficients. In addition to the harmonic coefficients, the values of the gravitational parameter $\mu$ and the equatorial radius of the central body $R_\oplus$ are what define a gravitational model \cite{WGS84,KonoplivBanerdtSjogren1999,KonoplivParkFolkner2016}. Due to the rotation of the central body, the gravitational potential depends on time when referred to an inertial frame. A customary simplification is that it rotates with constant rotation rate $n_\oplus$ about the polar axis, and hence the time dependence is avoided when the problem is formulated in a rotating frame.
\par

Reformulation of Eq.~(\ref{geopot}) in orbital elements is easily done by recalling the conversion from spherical to rectangular coordinates in the rotating frame
\begin{equation} \label{spherical}
\sin\phi=\frac{z}{r}, \qquad \sin\lambda=\frac{y}{\sqrt{x^2+y^2}}, \quad \cos\lambda=\frac{x}{\sqrt{x^2+y^2}},
\end{equation}
where $x$, $y$, and $z$, need to be expressed in the usual Keplerian elements: $a$, $e$, $I$, $\Omega$, $\omega$, $M$, for orbit semi-major axis, eccentricity, inclination, right ascension of the ascending node, argument of the periapsis, and mean anomaly, respectively. This is easily done by carrying out the rotations that relate the orbital frame with the Cartesian frame
\begin{equation} \label{rectangular}
\left(\begin{array}{c} x\\ y\\ z\end{array}\right)
=R_3(-h)\,R_1(-I)\,R_3(-\theta)\left(\begin{array}{c} r\\ 0\\ 0\end{array}\right),
\end{equation}
where $R_i$, $i=1,2,3$ denote standard rotation matrices, $h=\Omega-n_\oplus{t}$, is the argument of the ascending node of the orbit measured in the rotating frame, $\theta=f+\omega$ is the argument of the latitude, and $f$ is the true anomaly, which is an implicit function of $M$. Besides, the radius must be replaced in Eq.~(\ref{rectangular}) by the conic equation
\begin{equation} \label{conic}
r=\frac{p}{1+e\cos{f}},
\end{equation}
in which $p=a\eta^2$ is the parameter of the conic and $\eta=\sqrt{1-e^2}$ is customarily known as the eccentricity function.
\par

Direct substitution of Eqs.~(\ref{spherical})--(\ref{conic}) into Eq.~(\ref{geopot}) provides the required conversion of the gravitational potential into orbital elements, and is the usual procedure when dealing with just a few terms of the gravitational potential. However when the degree and order of the gravitational potential expansion go beyond the first few terms of Eq.~(\ref{geopot}), this transformation is notably expedited by using Kaula's developments \cite{Kaula1966}.
\par

On the other hand, the effects of tesseral perturbations are known to average out except for particular resonances of the satellite's mean motion with the rotation rate of the system. Then, we constrain to the zonal harmonics part of the gravity potential, which is obtained by neglecting in Eq.~(\ref{geopot}) those harmonic coefficients $S_{n,m}$, $C_{n,m}$ with $m>0$, viz.
\begin{equation} \label{zonalpot}
U=-\frac{\mu}{r}\sum_{n\ge{0}}\frac{R_\oplus^n}{r^n}C_{n,0}P_{n,0}(\sin\varphi).
\end{equation}
The zonal potential enjoys axial symmetry, and, therefore, is a 2-DOF problem that does not depend on the geocentric longitude. This reduction of the dimension releases the problem from the time dependency, and, therefore, the body-fixed frame is no longer needed and can be replaced by the usual inertial frame.
\par

Following Kaula's approach \cite{Kaula1966}, Eq.~(\ref{zonalpot}) is rewritten in orbital elements as
\begin{equation} \label{zonalKaula}
U=-\frac{\mu}{r}-\frac{\mu}{a}\left(\frac{a}{r}\right)^2\eta\sum_{i\ge2}V_i
\end{equation}
with
\begin{equation} \label{Vi}
V_i= \frac{R_\oplus^i}{a^i} \frac{{C}_{i,0}}{\eta^{2i-1}}
\sum_{j=0}^i\mathcal{F}_{i,j}(I)\sum_{k=0}^{i-1}\binom{i-1}{k}e^k\cos^k\!f \cos[(i-2 j)(f+\omega)-i_\pi],
\end{equation}
in which $\mathcal{F}_{i,j}$ are Kaula's inclination functions particularized for the case of the zonal problem, namely,
\begin{equation} \label{KaulaF}
\mathcal{F}_{i,j}=\sum_{l=0}^{\min(j,i_0)}\frac{(-1)^{j-l-i_0}}{2^{2i-2l}}\frac{(2i-2l)!}{l!(i-l)!(i-2l)!}\binom{i-2l}{j-l}\sin^{i-2l}I, \quad i\ge2l,
\end{equation}
and the parity correction $i_\pi$ and the symbol $i_0$ adhere to the index notation in \cite{Lara2018Stardust}. That is, for a generic summation index $i$
\begin{equation} \label{indexconvention}
i^\star=i\bmod2, \qquad 
i_\pi=\frac{\pi}{2}i^\star, \qquad 
i_m=\left\lfloor\frac{i-m}{2}\right\rfloor, \qquad 
i_m^\star=i_m+i^\star,
\end{equation}
where $m$ notes an integer number, and $\lfloor{p}/q\rfloor$ denotes the integer division of the integers $p$ and $q$. Besides, in what follows we shorten expressions by using the standard abbreviations $s\equiv\sin{I}$ and $c\equiv\cos{I}$.
\par

Note that Eqs.~(\ref{zonalKaula})--(\ref{Vi}) differ from analogous equations in \cite{Kaula1966} (see also \cite{RosboroughOcampo1991}) in which we keep a divisor $r^2$ factoring the summation comprising the non-centralities of the gravitational potential. The reason for that will be apparent later, in reference to the elimination of the parallax simplification.
\par

The flow derived from Eq.~(\ref{zonalKaula}) is obtained from the integration of Lagrange planetary equations \cite[see sect.~10.2 of][for instance]{Battin1999}, but analytical solutions to these equations are not generally known even for the lower degree truncations of Eq.~(\ref{zonalKaula}). However, under certain conditions, useful analytical approximations can be computed by perturbation methods \cite{Nayfeh2004,BrouwerClemence1961}.
\par

In particular, because the zonal problem is derived from a potential function, it admits Hamiltonian formulation and we can take benefit of Hamiltonian perturbations methods to approach the problem \cite{FerrazMello2007}. More specifically, we resort to the Lie transforms method as devised by Deprit \cite{Deprit1969}, which is summarized in  \ref{s:LieTransforms} for the convenience of interested readers. In Deprit's approach, the Hamiltonian of the zonal problem is written
\begin{equation} \label{zonalHam}
\mathcal{H}=\sum_{m\ge0}\frac{\epsilon^m}{m!}\mathcal{H}_{m,0},
\end{equation}
in which $\epsilon$ is a formal small parameter used to indicate the strength of each different perturbation, and we write the Hamiltonian terms in the form
\begin{eqnarray} \label{H00}
\mathcal{H}_{0,0} &=& -\frac{\mu}{2a}, \\ \label{H10}
\mathcal{H}_{1,0} &=& \left(-\frac{\mu}{a}\right)\frac{a^2\eta}{r^2}V_2, \\ \label{H20}
\mathcal{H}_{2,0} &=& 2\left(-\frac{\mu}{a}\right)\frac{a^2\eta}{r^2}\sum_{i\ge3}V_i, \\ \label{H30}
\mathcal{H}_{m,0} &=& 0, \qquad m\ge3,
\end{eqnarray}
\par

We note that orbital elements are not canonical variables, which are certainly required in Hamiltonian mechanics. However, because Keplerian elements provide an immediate insight into orbital problems, we regularly use the orbital elements notation in following expressions in the understanding that the symbols used are not variables, but known explicit or implicit functions of some set of canonical variables. In particular, the construction of perturbation solutions to perturbed Kepelrian problems is conveniently approached in Delaunay canonical variables $(\ell,g,h,L,G,H)$ corresponding to the mean anomaly $\ell=M$, the argument of the periapsis $g=\omega$, the argument of the node $h=\Omega$, the Delaunay action $L=\sqrt{\mu{a}}$, the total angular momentum $G=L\eta$, and the polar component of the angular momentum in the polar axis direction $H=G\cos{I}$ \cite{Delaunay1860}. 
\par

Note that the right ascension of the ascending node is a cyclic variable in the zonal Hamiltonian defined by the set of Eqs.~(\ref{Vi})--(\ref{H30}). For this reason, its conjugate momentum $H$ is an integral of the zonal problem. Common perturbation solutions to the zonal problem are based in finding a canonical transformation that completely removes the mean anomaly $\ell$ from the Hamiltonian. After truncation to the order of the transformation, the reduced Hamiltonian is of 1-DOF, and, in consequence, it is integrable in the prime variables.

However, one must note that, because of the peculiarities of the elliptic motion the mean anomaly does not appear explicitly in the zonal Hamiltonian, but implicitly through the true anomaly $f\equiv{f}(\ell,G,L)$. This fact makes the equation of the center $\phi=f-\ell$ to appear in the computation of the generating function as soon as in the first order of a perturbation approach, in this way making notably difficult the computation of higher orders of the solution in closed form \cite{Kozai1962,OsacarPalacian1994}. The traditional alternative relies on expansions of the elliptic motion using trigonometric functions \cite{BrouwerClemence1961}, but this approach requires high order expansions for effectively dealing with moderate eccentricities, with the consequent long series to evaluate \cite{DepritRom1970,Wnuk1997}, and does not solve efficiently the case of high eccentricities. Even though alternative expansions of the elliptic motion have been suggested using elliptic function theory \cite{BrumbergFukushima1994} leading to extensions of Kaula's recursions \cite{Brumbergetal1995}, they are rarely used due to the technicalities involved in the evaluation of the required special functions \cite{Fukushima2014ch}. Therefore, the standard way in the reduction of the zonal Hamiltonian starts from making a preparatory simplification which is traditionally dubbed \emph{the elimination of the parallax}  \cite{Deprit1981,LaraSanJuanLopezOchoa2013b,LaraSanJuanLopezOchoa2013c}.

\section{Simplification: elimination of the parallax} \label{s:parallax}

The aim of this simplification is to remove non-essential short-period terms from the zonal Hamiltonian in Eq.~(\ref{zonalHam}). That is, with this simplification, all the short-period terms are removed except those which are related to the fundamental Kepler equation; hence the flow derived from the in this way simplified Hamiltonian is termed quasi-Keplerian motion \cite{Deprit1981}.
\par

The elimination of the parallax simplification applies to full zonal models \cite{AlfriendCoffey1984}, but we don't need to do that due to the low order of the perturbation theory we are interested in. Hence, we limit the application of this simplification just to case of the zonal harmonic of the second degree. That is, short period terms related with the explicit appearance of $f$ in Eq.~(\ref{H10}), which happen in the $V_2$ coefficient, as shown in Eq.~(\ref{Vi}), will be removed by the transformation up to the second order of $C_{2,0}$, whereas the explicit appearance of $r$ in Eq.~(\ref{H10}) , which is also affected of short-period terms, cf.~Eq.~(\ref{conic}), will remain untouched by the transformation, as will remain so the remaining terms of the Hamiltonian, which are given in Eq.~(\ref{H20}).

\subsection{1st order}

Then, the first step in solving the homological equation of Deprit's perturbation approach in \ref{s:LieTransforms}, is to choose the new Hamiltonian term $\mathcal{H}_{0,1}$ to be comprised of those terms of $\mathcal{H}_{1,0}$ in Eq.~(\ref{H10}) that are free from the explicit appearance of $f$. In view of
\begin{equation}
V_2= -\frac{R_\oplus^2}{a^2}\frac{C_{2,0}}{\eta^3}\Big\{\frac{1}{4}\left(2-3s^2\right) (1+e\cos{f})
+\frac{3}{8}s^2\left[e\cos(f+2\omega)+2\cos(2f+2\omega)+e\cos(3f+2\omega)\right] \Big\},
\end{equation}
we choose
\begin{equation} \label{H01parallax}
\mathcal{H}_{0,1}=-\frac{\mu}{a}\frac{a^2\eta}{r^2}\langle{V}_2\rangle_f,
\end{equation}
with
\begin{equation} \label{V2averaged}
\langle{V}_2\rangle_f = -C_{2,0}\frac{R_\oplus^2}{p^2}\frac{1}{4}\eta\left(2-3s^2\right).
\end{equation}
\par

Then, $W_1$ is computed from Eq.~(\ref{homoDelo}) by simple quadrature as
\[
W_1=\frac{1}{n}\int(\mathcal{H}_{1,0}-\mathcal{H}_{0,1})\,\text{d}\ell, 
\]
that can be solved in closed form using the differential relation between the true and mean anomalies
\begin{equation} \label{dldf}
a^2\eta\,\text{d}\ell=r^2\,\mathrm{d}f,
\end{equation}
which is derived from the preservation of the angular momentum in the Keplerian motion. Thus,
\[
W_1=\frac{1}{n}\int(\mathcal{H}_{1,0}-\mathcal{H}_{0,1})\,\frac{r^2}{a^2\eta}\,\text{d}f
=-L\int\left(V_2-\langle{V}_2\rangle_f\right)\,\text{d}f,
\]
and hence
\begin{equation} \label{W1parallax}
W_1=G\frac{R_\oplus^2}{p^2}\frac{C_{2,0}}{8}\left[e\left(4-6 s^2\right)\sin{f}+s^2\sum_{i=0}^3Q_i(e)\sin (if+2g)\right],
\end{equation}
with
\begin{equation} \label{Q0123}
Q_0
=\frac{1+2\eta}{(1+\eta)^2}e^2, \qquad Q_1=3e, \qquad Q_2=3, \qquad Q_3=e.
\end{equation}
\par

Note that the summand of $W_1$ related to the term $Q_0$ does not depend on $\ell$. It is an integration ``constant'' that has been added to the integral to guarantee that $W_1$ is free from hidden long-period terms on the argument of the periapsis \cite{Kozai1962,ExertierThesis1988}.\footnote{These kinds of functions that are free from the mean anomaly are sometimes called ``Kozai-like constants'' in artificial satellite theory \cite{LaraSanJuanLopezCefola2012}.} Indeed, it is simple to check that $\int_{0}^{2\pi}W_1\,\mathrm{d}\ell=0$ in this way being free from long-period terms. With the addition of this integration constant to the first order of the generating function we guarantee that the second order Hamiltonian, still to be computed, will not be deprived of any long-period effect \cite{LaraSanJuanLopezOchoa2013b,LaraSanJuanLopezOchoa2013c}. 
\par

We hasten to add that dealing with integration constants of the generating function is not necessary from the point of view of constructing the perturbation solution, which involves the transformation from original to new variables and vice versa, yet it is highly recommended in those cases in which realistic information is expected to be obtained directly from the simplified Hamiltonian.

\subsection{2nd order}

At the second order, Eq.~(\ref{homoDelo}) leads to
\begin{equation} \label{homo2parallax}
W_2=\frac{1}{n}\int(\widetilde{\mathcal{H}}_{0,2}-\mathcal{H}_{0,2})\,\mathrm{d}\ell,
\end{equation}
where, from Deprit's recursion in Eq.~(\ref{deprittriangle}),
\begin{equation} \label{H02tilde}
\widetilde{\mathcal{H}}_{0,2}=\{\mathcal{H}_{0,1},W_1\}+\{\mathcal{H}_{1,0},W_1\}+\mathcal{H}_{2,0}.
\end{equation}
After evaluating the Poisson brackets and in view of Eq.~(\ref{H20}), it is obtained
\begin{equation} \label{tildeH02parallax}
\widetilde{\mathcal{H}}_{0,2} = \frac{\mu}{a}\frac{a^2\eta}{r^2}\left\{ 
C_{2,0}^2\frac{R_\oplus^4}{p^4}\eta
\sum_{j=0}^6\sum_{l=-2}^2\left[q_{j,l}+\frac{es^2}{(1+\eta)^2}\tilde{q}_{j,l}\right]\cos(jf+2l\omega) 
-2\sum_{i\ge3}V_i
\right\}
\end{equation}
in which the non-null coefficients $q_{j,l}\equiv{q}_{j,l}(s,e,\eta)$ and $\tilde{q}_{j,l}\equiv\tilde{q}_{j,l}(s,e,\eta)$ are presented in Tables \ref{t:qjlparallax} and \ref{t:tildeqjlparallax}, respectively.\footnote{For computational purposes, one may note that, following the index convention in Eq.~(\ref{indexconvention}), the lower limit of the summation index $l$ for the $q_{j,l}$ coefficients is $l=j_1$, whereas in the case of $\tilde{q}_{j,l}$ is $l=j_4$.} Recall that the term $\tilde{q}_{0,\pm1}$ is a consequence of including the term $Q_0$ of Eq.~(\ref{Q0123}) in the first order generating function, and can be neglected when having ``centered'' mean elements is not a requirement of the perturbation theory.
\par
\begin{table*}[htb] %
\centering 
\begin{tabular}{@{}llll@{}}
$j$ & \multicolumn{1}{l}{$l=0$} & \multicolumn{1}{l}{$l=1$} & \multicolumn{1}{l}{$l=2$}  \\[0.33ex]
\hline
$0$ & $q_{2,0}-\frac{1}{16}(21 s^4-42 s^2+20)$ & $\frac{3}{64} e^2 s^2 \left(14-15 s^2\right)$ & \vphantom{$\frac{M^9}{M^9}$} 
\\ [0.5ex]
$1$ & $-\frac{1}{32}e\left(27s^4-108s^2+64\right)$ & $\frac{7}{16}es^2\left(11-12s^2\right)$ &  
\\ [0.5ex]
$2$ & $\frac{3}{64}e^2\left(5s^4+8s^2-8\right)$ & $\frac{3}{16}e^2s^2\left(2-s^2\right)+\frac{1}{8}\left(20-21s^2\right)s^2$ & $-\frac{15}{128}e^2s^4$ 
\\ [0.5ex]
$3$ & & $\frac{3}{16}es^2\left(8s^2-5\right)$ & $-\frac{9}{64}es^4$ 
\\ [0.5ex]
$4$ & & $\frac{3}{32}e^2s^2\left(13s^2-10\right)$ & $\frac{3}{64}\left(4-e^2\right)s^4$ 
\\ [0.5ex]
$5$ & & & $\phantom{-}\frac{15}{64}es^4$ 
\\ [0.5ex]
$6$ & & & $\phantom{-}\frac{9}{128}e^2s^4$ 
\\ [0.33ex]
\hline
\end{tabular}
\normalsize
\caption{Coefficients $q_{j,l}$ in Eqs.~(\protect\ref{tildeH02parallax}).
}
\label{t:qjlparallax}
\end{table*}
\begin{table*}[htb] %
\centering 
\begin{tabular}{@{}rlll@{}}
$l$ & \multicolumn{1}{l}{$j=1$} & \multicolumn{1}{l}{$j=2$}  \\[0.33ex]
\hline
$-2$ & $-\frac{3}{256}e^2s^2\chi^{-}$ &  & \vphantom{$\frac{M^9}{M^9}$} 
\\ [0.75ex]
$-1$ & $\frac{3}{128}\left[(9-42\eta-31\eta^2)s^2-\frac{2}{3}(5-54\eta-39\eta^2)\right]\chi^{-}$
     & $\frac{1}{8}e(3s^2-2)\chi^{-}$
\\ [0.75ex]
 $0$ & $\frac{3}{128}(2+23\eta-31\eta^2)s^2\chi^{+}-\frac{3}{16}e^2(1+2\eta)$
     & $\frac{3}{32} e[(15s^2-8)\eta-4]$
\\ [0.75ex]
$+1$ & $\frac{3}{128}[(37\eta^2-14\eta-83)s^2+(58+12\eta-30\eta^2)]\chi^{+}$
     & $-\frac{3}{8}e(3s^2-2)\chi^{+}$
\\ [0.75ex]
$+2$ & $\frac{3}{256}(21+18\eta+\eta^2)s^2\chi^{-}$
     & $\frac{15}{32} e (\eta +2) s^2$
\\ [0.33ex]
\hline \\[-1ex]
$l$ & \multicolumn{1}{l}{$j=3$} & \multicolumn{1}{l}{$j=4$} & \multicolumn{1}{l}{$j=5$}  \\[0.33ex]
\hline
$-1$ & $\frac{3}{16}e\tilde{q}_{2,-1}$  & \vphantom{$\frac{M^9}{M^9}$} 
\\ [0.75ex]
 $0$ & $\frac{3}{256}(61\eta^2+66\eta-23)s^2\chi^{-}-\frac{3}{16}e^2(1+2\eta)$ 
     & $-\frac{9}{32}es^2\chi^{-}$
     & $\frac{5}{24}e\tilde{q}_{4,0}$
\\ [0.75ex]
$+1$ & $\frac{3}{16}e\tilde{q}_{2,1}$ 
\\ [0.75ex]
$+2$ & $\frac{9}{256}(39-6\eta-5\eta^2)s^2\chi^{+}$ 
     & $\frac{27}{32}es^2\chi^{+}$
     & $\frac{5}{24}e\tilde{q}_{4,2}$
\\ [0.33ex]
\hline
\end{tabular}
\normalsize
\caption{Coefficients $\tilde{q}_{j,l}$ in Eqs.~(\protect\ref{tildeH02parallax}); $\tilde{q}_{0,\pm1}=\frac{3}{16}e(1+2\eta)(4-5s^2)$ and $\chi^{\pm}\equiv1\pm\eta$.}
\label{t:tildeqjlparallax}
\end{table*}

In the same way as in the first order, the new Hamiltonian term $\mathcal{H}_{0,2}$ is chosen to be made of those terms of $\widetilde{\mathcal{H}}_{0,2}$ that are free from the explicit appearance of $f$ but keep the explicit appearance of $r$ in Eq.~(\ref{tildeH02parallax}) untouched. Hence 
\begin{equation} \label{H02parallax}
\mathcal{H}_{0,2} = 2\frac{\mu}{a}\frac{a^2\eta}{r^2}\left\{
C_{2,0}^2\frac{R_\oplus^4}{p^4}\eta
\left[\frac{q_{0,0}}{2}+\left(q_{0,1}+\frac{es^2\tilde{q}_{0,1}}{(1+\eta)^2}\right)\!\cos2\omega\right]
-\sum_{i\ge3}\langle{V}_i\rangle_f\right\},
\end{equation}
where $\langle{V}_i\rangle_f$ comprises the terms of Eq.~(\ref{Vi}) that are free from short-period effects depending on $f$. It is computed as follows.
\par

The dependency of $V_i$ on the true anomaly is made explicit by expanding Eq.~(\ref{Vi}) as a Fourier series in $f$. This is done by expanding first
\begin{equation*}
\cos^kf\cos[m(f+\omega)-i_\pi]=
\cos^kf\cos{mf}\cos(m\omega-i_\pi)-\cos^kf\sin{mf}\sin(m\omega-i_\pi),
\qquad m=(i-2j).
\end{equation*}
Then, using the standard trigonometric reductions
\begin{equation*}
\cos\beta\cos^k\alpha=\frac{\cos (\alpha +\beta )+\cos (\alpha -\beta )}{2}\cos^{k-1}\alpha, \qquad
\sin\beta\cos^k\alpha=\frac{\sin (\alpha +\beta )-\sin (\alpha -\beta )}{2}\cos^{k-1}\alpha,
\end{equation*}
one easily arrives to the recursion
\begin{equation} \label{recurre}
\cos^k\!f\left\{\begin{array}{l} \cos(i-2 j)f \\ \sin(i-2 j)f \end{array}\right\} = 
\frac{1}{2^k}\sum_{l=0}^k\binom{k}{k-l}\left\{\begin{array}{l} \cos(i-2j-k+2l)f \\ \sin(i-2j-k+2l)f \end{array}\right\},
\end{equation}
which shows that the only terms in Eq.~(\ref{Vi}) that are free from $f$ come from those terms of Eq.~(\ref{recurre}) such that $l=\frac{1}{2}[k-(i-2j)]$ or $k-l=\frac{1}{2}(k+i)-j$. Hence,
\begin{equation} \label{VinofG}
\langle{V}_i\rangle_f= \frac{R_\oplus^i}{a^i}C_{i,0}\sum_{j=0}^i\mathcal{F}_{i,j}(s)\mathcal{G}_{i,j}(e,\eta)\cos[(i-2j)\omega-i_\pi],
\end{equation}
where
\begin{equation} \label{miG}
\mathcal{G}_{i,j}=\frac{1}{\eta^{2i-1}} \sum_{k=0}^{i-1}\binom{i-1}{k}\binom{k}{\frac{k+i}{2}-j}\frac{e^k}{2^k}.
\end{equation}
\par

As expected, the functions $\mathcal{G}_{i,j}$ in Eq.~(\ref{miG}) are no more than Kaula's eccentricity functions for the particular case of the zonal problem. Indeed, because the summation index $l$ in Eq.~(\ref{recurre}) is integer, $k+i$ must be even in Eq.~(\ref{miG}), which, therefore, can be rearranged in the efficient form
\begin{equation} \label{KaulaG}
\mathcal{G}_{i,j}=\frac{1}{\eta^{2i-1} }\sum_{l=0}^{\tilde\jmath-1}\binom{i-1}{q}\,\binom{q}{l}\,\frac{e^q}{2^q},
\qquad q=2l+i-2\tilde\jmath, \qquad 
\left\{ \begin{array}{rcl} i\ge2j & \Rightarrow & \tilde\jmath=j \\[0.5ex] i<2j & \Rightarrow & \tilde\jmath=i-j \end{array}, \right. 
\end{equation}
proposed by Kaula, cf.~Eq. 3.66 of \cite{Kaula1966}.
\par

Taking into account the peculiarities of the summations for the cases of odd and even harmonics, Eq.~(\ref{VinofG}) is more efficiently organized in the form
\begin{equation} \label{Vistar}
\langle{V}_i\rangle_f=
\frac{R_\oplus^i}{a^i}C_{i,0}\sum_{j=0}^{i_2}
(2-\delta_{j+i^\star,0})\mathcal{F}_{i,i_0^\star+j}(s)\mathcal{G}_{i,i_0^\star+j}(e,\eta)\cos\left[(2j+i^\star)\omega+i_\pi\right],
\end{equation}
in which $\delta_{i,j}$ is the usual Kronecker delta function. Note that Eq.~(\ref{Vistar}) also applies to the term $\langle{V}_2\rangle_f$, which was displayed expanded in Eq.~(\ref{V2averaged}).
\par

\section{Delaunay normalization} \label{s:delaunay}

The transformed, simplified Hamiltonian after the elimination of the parallax is obtained by replacing the original Delaunay variables by new, prime ones, viz. $(\ell',g',h',L',G',H')$. Hence, as follows from Eqs.~(\ref{H01parallax}) and (\ref{H02parallax}), up to the second order, the new Hamiltonian is
\[
\mathcal{H}'=\mathcal{H}'_{0,0}+\epsilon\mathcal{H}'_{1,0}+\frac{\epsilon^2}{2}\mathcal{H}'_{2,0}+\mathcal{O}(C_{2,0}^3)
\]
with
\begin{eqnarray} \label{H00p}
\mathcal{H}'_{0,0} &=& -\frac{\mu}{2a}, \\ \label{H10p}
\mathcal{H}'_{1,0} &=& -\frac{\mu}{a}\frac{a^2\eta}{r^2}\langle{V}_2\rangle_f, \\  \label{H20p}
\mathcal{H}'_{2,0} &=& 2\frac{\mu}{a}\frac{a^2\eta}{r^2}\left\{C_{2,0}^2\frac{R_\oplus^4}{p^4}\eta
\left[\frac{q_{0,0}}{2}+\left(q_{0,1}+\frac{es^2}{(1+\eta)^2}\tilde{q}_{0,1}\right)\cos2\omega\right]
-\sum_{i\ge3}\langle{V}_i\rangle_f\right\},
\end{eqnarray}
where $\langle{V}_2\rangle_f$ and $\langle{V}_i\rangle_f$ are given in Eqs.~(\ref{V2averaged}) and (\ref{Vistar}), respectively, and all symbols are now assumed to be functions of the prime Delaunay variables.
\par

The elimination of the mean anomaly now becomes trivial. Indeed, at the first order we choose
\begin{equation} \label{H01Delaunay}
\mathcal{H}_{0,1}'=\frac{1}{2\pi}\int_0^{2\pi}\mathcal{H}'_{1,0}\,\mathrm{d}\ell'
=\frac{1}{2\pi}\int_0^{2\pi}\mathcal{H}'_{1,0}\frac{r^2}{a^2\eta}\,\mathrm{d}f
=-\frac{\mu}{a}\langle{V}_2\rangle_f.
\end{equation}
The first order term of the generating function is computed from Eq.~(\ref{homoDelo})
\[
W'_1=\frac{1}{n}\int(\mathcal{H}'_{1,0}-\mathcal{H}'_{0,1})\,\text{d}\ell'
=\frac{1}{n}\left[-\mathcal{H}'_{0,1}\ell'+\int\mathcal{H}'_{1,0}\frac{r^2}{a^2\eta}\,\text{d}f\right]
=-L\langle{V}_2\rangle_f\phi
\]
where $\phi=f-\ell'$ is the equation of the center. 
\par

The term $\widetilde{\mathcal{H}}'_{0,2}$ is computed again from Deprit's recursion in Eq.~(\ref{deprittriangle}), which gives formally the same expression as Eq.~(\ref{H02tilde}), but now it must be formulated in the tilde functions and variables. The new Hamiltonian term $\mathcal{H}'_{0,2}$ is selected to be free from short-period terms, viz.
\[
\mathcal{H}'_{0,2}=\frac{1}{2\pi}\int_0^{2\pi}\widetilde{\mathcal{H}}'_{0,2}\,\mathrm{d}\ell'
=\mathcal{H}'_{0,2,1}+\mathcal{H}'_{0,2,2}+\mathcal{H}'_{0,2,3},
\]
where
\begin{eqnarray*}
\mathcal{H}'_{0,2,1} &=& \frac{1}{2\pi}\int_0^{2\pi}\left(\{\mathcal{H}'_{0,1},W'_1\}+\{\mathcal{H}'_{1,0},W'_1\}\right)\text{d}\ell'
\\
\mathcal{H}'_{0,2,2} &=& \frac{1}{2\pi}\int_0^{2\pi}
2\frac{\mu}{a}\frac{a^2\eta}{r^2}C_{2,0}^2\frac{R_\oplus^4}{p^4}\eta
\left\{\frac{q_{0,0}}2+\left[q_{0,1}+\frac{es^2}{(1+\eta)^2}\tilde{q}_{0,1}\right]\!\cos2\omega\right\}\,\text{d}\ell',
\\
\mathcal{H}'_{0,2,3} &=& \frac{1}{2\pi}\int_0^{2\pi}
\Bigg(-2\frac{\mu}{a}\frac{a^2\eta}{r^2}\sum_{i\ge3}\langle{V}_i\rangle_f
\Bigg)\,\text{d}\ell',
\end{eqnarray*}
After evaluating the Poisson brackets, and using, when required, the differential relation in Eq.~(\ref{dldf}), the quadratures are solved in closed form of the eccentricity to get
\begin{eqnarray*}
\mathcal{H}'_{0,2,1} &=& -\frac{\mu}{2a}C_{2,0}^2\frac{R_\oplus^4}{p^4}\frac{1}{8}\eta(1+3\eta)(2-3s^2)^2,
\\
\mathcal{H}'_{0,2,2} &=& 2\frac{\mu}{a}C_{2,0}^2\frac{R_\oplus^4}{p^4}\eta\left\{
\frac{q_{0,0}}{2}+\left[q_{0,1}+\frac{es^2}{(1+\eta)^2}\tilde{q}_{0,1}\right]\!\cos2\omega\right\},
\\
\mathcal{H}'_{0,2,3} &=& -2\frac{\mu}{a}\sum_{i\ge3}\langle{V}_i\rangle_f.
\end{eqnarray*}
Hence
\begin{eqnarray} \label{H02Delaunay}
\mathcal{H}'_{0,2} &=& -\frac{\mu}{a}2\sum_{i\ge3}\langle{V}_i\rangle_f
-\frac{\mu}{a}\frac{R_\oplus^4}{p^4}C_{2,0}^2\eta
\left\{\frac{1+3\eta}{16}(2-3s^2)^2-q_{0,0}-2\left[q_{0,1}+\frac{es^2\tilde{q}_{0,1}}{(1+\eta)^2}\right]\cos2\omega\right\}.
\end{eqnarray}
\par

Up to $\mathcal{O}(\epsilon^3)=\mathcal{O}(C_{2,0}^3)$, the new Hamiltonian is
\begin{equation} \label{Hampp}
\mathcal{H}''=\mathcal{H}''_{0,0}+\epsilon\mathcal{H}''_{1,0}+\frac{\epsilon^2}{2}\mathcal{H}''_{2,0},
\end{equation}
where terms $\mathcal{H}''_{m,0}$ $(m=0,1,2)$, are obtained from corresponding ones $\mathcal{H}'_{0,m}$ given by Eqs.~(\ref{H00p}), (\ref{H01Delaunay}), and (\ref{H02Delaunay}), respectively, by changing prime by double prime variables $(\ell'',g'',h'',L'',G'',H'')$.
\par

Hence, if we neglect the term $\tilde{q}_{0,1}$ from Eq.~(\ref{H02Delaunay}) ---which is a consequence of the integration constant introduced in Eq.~(\ref{W1parallax}) at the 1st order of the elimination of the parallax simplification--- the expanded Hamiltonian in Appendix A of \cite{CoffeyDepritDeprit1994} can be replaced by the compact expression
\[
\mathcal{H}=
\frac{G^2}{p^2}\left\{C_{2,0}^2\eta^3\frac{R_\oplus^4}{p^4}
\left[-\frac{1+3\eta}{32}(2-3s^2)^2+\frac{1}{2}q_{0,0}+q_{0,1}\cos2\omega\right]
-\eta^2\sum_{i\ge2}\langle{V}_i\rangle_f\right\},
\]
together with the use of Kaula's recursion in Eq.~(\ref{Vistar}).

\section{Kaula-type recursion's performance} \label{s:performance}

The Hamiltonian scaling used in Eqs.~(\ref{H00})--(\ref{H30}) is valid only in the case of earth-like bodies \cite{CoffeyDepritDeprit1994}. Interesting as this case may be, there are cases, as it happens with Venus or the moon gravitational potentials, in which the prevalence of the $C_{2,0}$ coefficient is not enough to define a clear scaling of the harmonic coefficients, which, therefore, are all taken to be of the same order in the perturbation arrangement. In these cases, 2nd order effects of the second zonal harmonic are not needed and the elimination of the parallax simplification can be avoided. Then, the mean elements Hamiltonian (\ref{Hampp}) is obtained after a single canonical transformation that results in
\begin{equation} \label{Hmean}
\mathcal{H}'=- \frac{\mu}{2a}-\frac{\mu}{a}\sum_{i\ge2}\langle{V}_i\rangle_f,
\end{equation}
where $\langle{V}_i\rangle_f$ is given in Eq.~(\ref{Vistar}). The derivation and arrangement of Eqs.~(\ref{Hmean}) is definitely simpler than corresponding expressions in \cite{Saedeleer2005}.
\par

\begin{figure}[htb]
\centering
\includegraphics[scale=0.7]{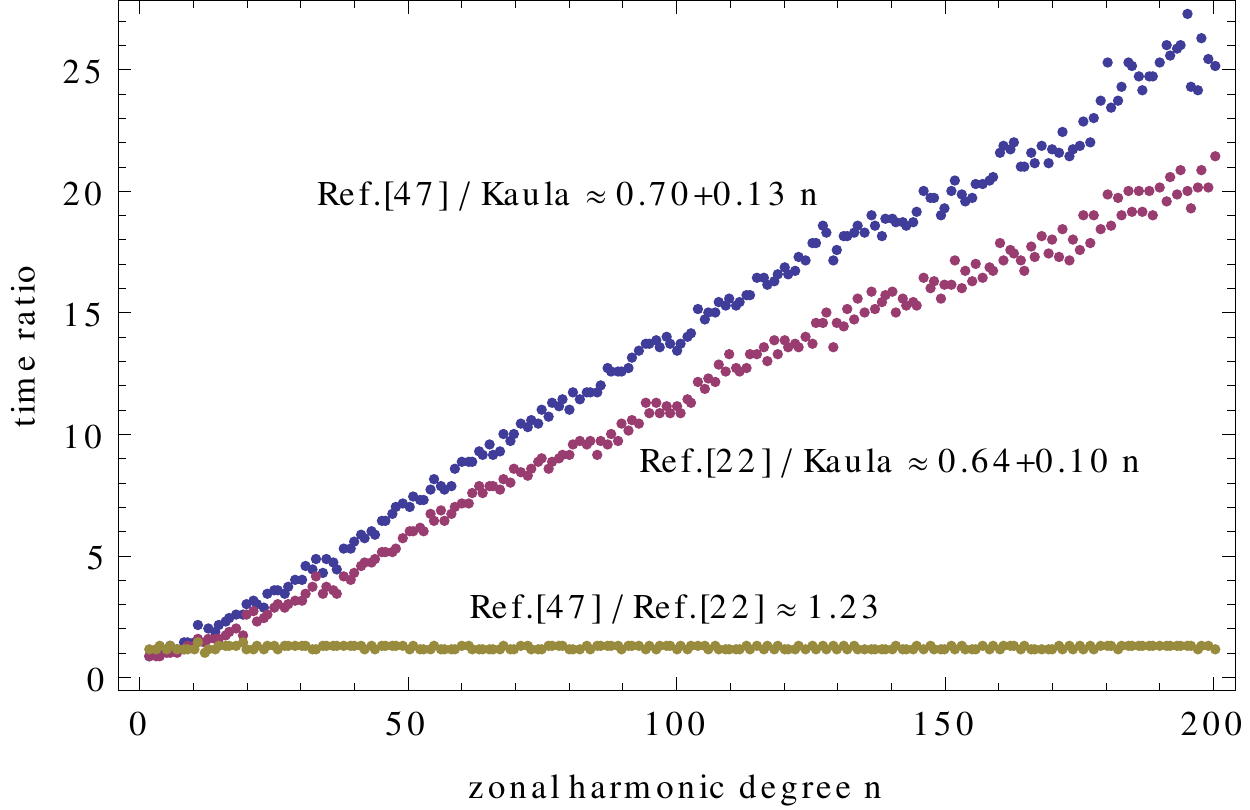}
\caption{Performance of Eq.~(\ref{Hmean}) in the construction of the mean elements Hamiltonian compared to Eqs.~(49)--(50) of \cite{Saedeleer2005} and corresponding improved recursions in \cite{LaraSaedeleerFerrer2009}.
}
\label{f:Ham}
\end{figure}

We compare the efficiency of Eq.~(\ref{Hmean}) in the construction of the mean elements Hamiltonian with corresponding ones proposed in \cite{Saedeleer2005}, but also with the improved version of these equations developed in \cite{LaraSaedeleerFerrer2009}. The results of the comparisons are summarized in Fig.~\ref{f:Ham}, where the ratio of the time needed in the construction of a given term $n$ of the averaged potential is displayed for the three following cases:
\begin{itemize}
\item Eqs.~(48)--(50) of Ref.~\cite[][p.~250]{Saedeleer2005} compared to Eq.~(\ref{Hmean}), represented by blue dots,
\item improved equations in Ref.~\cite{LaraSaedeleerFerrer2009} compared to Eq.~(\ref{Hmean}), represented by magenta dots,
\item Eqs.~(48)--(50) of Ref.~\cite{Saedeleer2005} compared to improved equations in Ref.~\cite{LaraSaedeleerFerrer2009}, represented by yellow dots.
\end{itemize}
As shown in the figure, all the approaches spend similar times for the lower degrees of the potential. Besides, the recursions in \cite{LaraSaedeleerFerrer2009} reduce the time needed in the computation of the original recursions in \cite{Saedeleer2005} to about the 80\%. But even in that case, De Saedeleer's recursions spend a clearly increasing time, with a rate that is approximately proportional to 1 tenth of the zonal harmonic degree $n$, with respect to Kaula's recursions. Thus, in De Saedeleer's (improved) approach, the time spent in constructing the zonal term of 50th degree is $\sim5$ times longer than the time needed when using Kaula's recursions, $\sim10$ times longer for the zonal term of 100th degree, or $\sim20$ times longer for the zonal term of 200th degree.
\par

\begin{figure}[htb]
\centering
\includegraphics[scale=0.7]{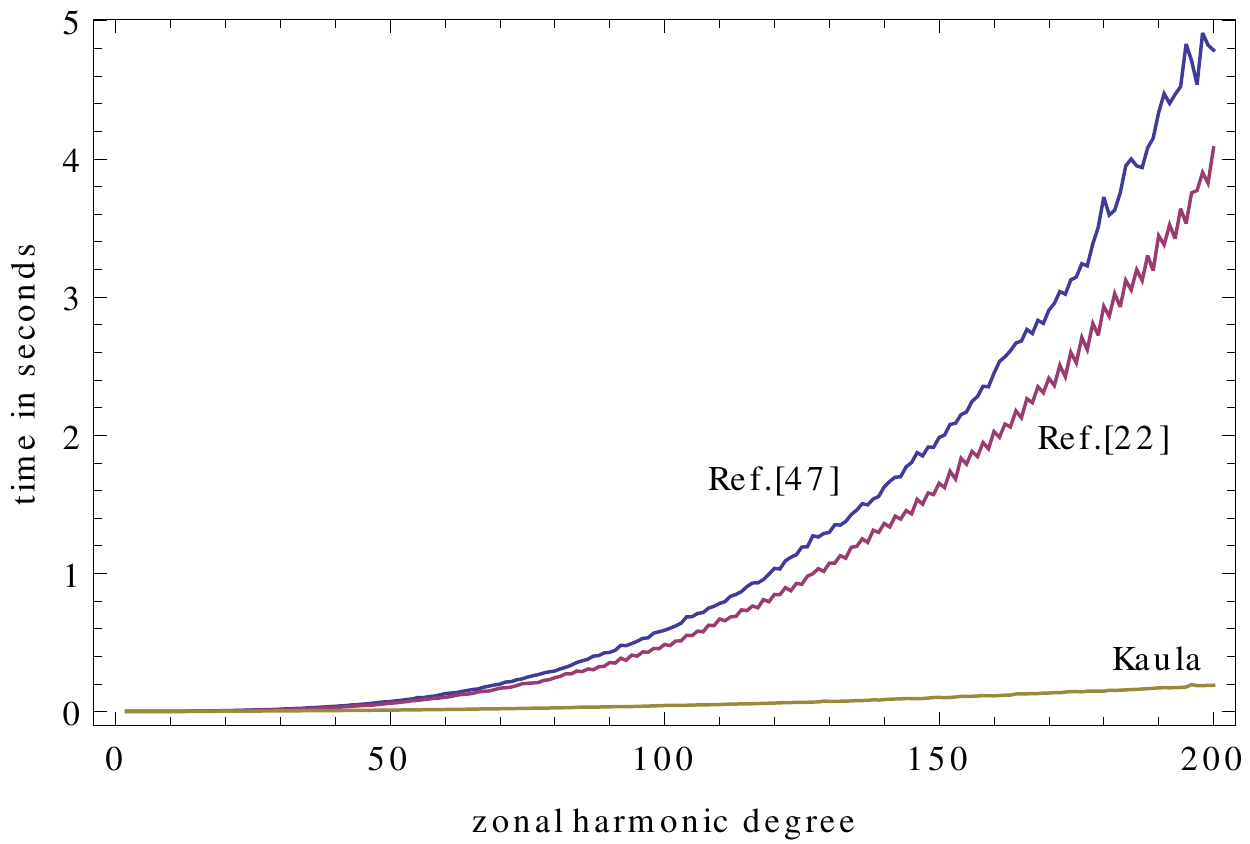}
\caption{Time spent in the computation of each term of the mean elements Hamiltonian with the different recursions.
}
\label{f:Htimes}
\end{figure}

In spite of the clear advantage of using Kaula's recursions in the construction of the mean elements Hamiltonian, all the approaches are quite feasible with the current computational power. Indeed, as shown in Fig.~\ref{f:Htimes}, even for the higher degrees, the analytical expression of the mean elements Hamiltonian term is obtained with Wolfram \textit{Mathematica}\raisebox{.9ex}{\tiny\textregistered} 9 software in just a few seconds when using a 2.8 GHz Intel Core i7 with 16 GB of RAM running under macOS High Sierra 10.13.5. However, while computing time grows with the degree in a cubic rate when using De Saedeleer's approach, it just grows only slightly higher than quadratic when using Kaula's recursions.

\section{The brut force approach}

When the Hamiltonian normalization process is approached directly by perturbations with the help of a computer algebra system, the highest efficiency is obtained when all the analytical expression to be handled in the computation of the Poisson brackets, as well as other operators that feed the Lie transform procedure, are fully expanded. Because of that, the number of terms required by particular perturbation theories (Hamiltonian, generating function, periodic corrections) are repeatedly reported in the literature \cite[see][among others]{DepritRom1970,CoffeyAlfriend1981,CoffeyDeprit1982,CoffeyNealSegermanTravisano1995,LaraSanJuanLopezOchoa2013c}. Once the perturbation problem has been solved, there is no doubt that the perturbation solution must be exactly the same as the one obtained with the recursions. However, the latter provides the solution in a more organized way, in which the necessary expressions take the structure of Poisson series \cite{DanbyDepritRom1965,DepritRom1968} where the coefficients of the trigonometric functions are arranged in the form of eccentricity and inclination polynomials, contrary to expanded monomials, as can be checked in Eq.~(\ref{Vistar}).
\par

In spite of the Lie transforms procedure is straightforward in the case of a simple first order approach, which was precisely the case tested in Section \ref{s:performance}, the fact that it starts from reformulating the zonal gravitational potential in orbital elements, which is given in Eq.~(\ref{zonalpot}) in spherical coordinates, hampers the brut force approach since the beginning. The whole procedure is relieved from unnecessary computations when starting from Eqs.~(\ref{zonalKaula})--(\ref{Vi}), which notably ease the computation of the mean elements Hamiltonian in closed form of the eccentricity using the differential relation in Eq.~(\ref{dldf}). But even when this shortcut is taken, the computations involved in the customary reduction of the powers of trigonometric functions that appear in Eq.~(\ref{Vi}) to trigonometric functions with combined arguments, very soon need to handle notable amounts of memory, in this way making the computation time to grow in a quintic way from degree to degree.
\par

\begin{figure}[htb]
\centering
\includegraphics[scale=0.7]{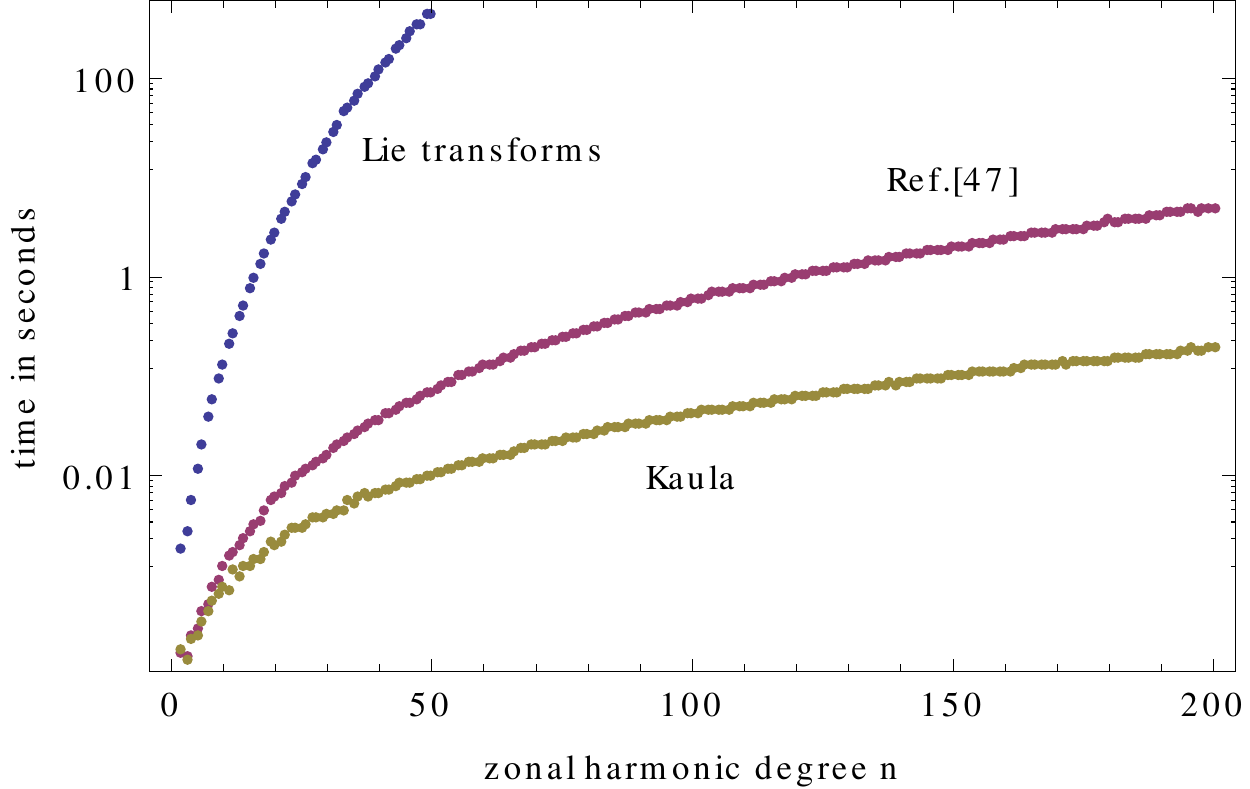}
\caption{Time spent in the computation of each term of the mean elements Hamiltonian with the brut force approach (Lie transforms), and De Saedeleer's and Kaula's recursions.
}
\label{f:LieHamt}
\end{figure}

\begin{figure}[htb]
\centering
\includegraphics[scale=0.7]{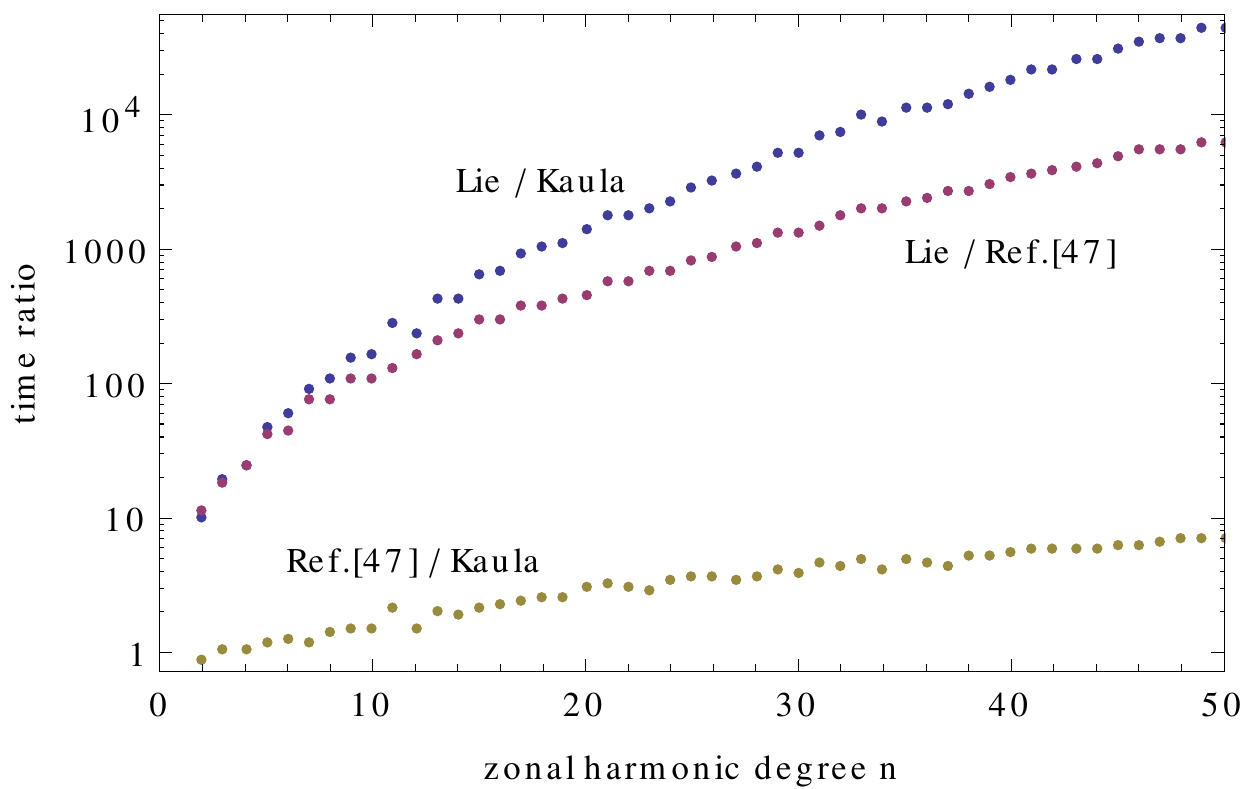}
\caption{Performance of Eq.~(\ref{Hmean}) in the construction of the mean elements Hamiltonian compared to Eqs.~(49)--(50) of \cite{Saedeleer2005} and corresponding improved recursions in \cite{LaraSaedeleerFerrer2009}.
}
\label{f:HamLie}
\end{figure}

The close to exponential growth of the Lie transforms approach is illustrated in Fig.~\ref{f:LieHamt}, in a logarithmic scale, where the computations with the Lie transforms where extended only to the 50th degree term of the zonal gravitational potential. The Lie transforms performance relative to both types of recursions is depicted in Fig.~\ref{f:HamLie}, where the performance of Kaula's recursions over those in \cite{Saedeleer2005} is also shown as reference. It can be observed in Fig.~\ref{f:HamLie} that, since the beginning, the performance is at least 10 times better on the side of recursions, it immediately grows to 100 times better even for the lower degrees, and very soon grows to 1000 times better when the degree is only moderate.


The numbers obtained in our comparisons with the direct implementation of the perturbation method by Lie transforms are just provided as an illustration of the performance, because the perturbation approach may admit different implementations. We used our long experience with the Lie transforms method to make the comparisons unbalanced, yet in no way we claim that the procedure cannot be speeded. However, due to the results obtained, it is not expected that more efficient implementations of the perturbation approach, if possible, will notably change our reported results. 
\par

Also, it may be argued against our comparisons that we relied on a general purpose algebra system whereas the computations can be done much more efficiently if this task is programmed in a specific symbolic manipulator. This is obviously true, but analogous improvements would be found in the recursions if they are likewise programmed with a specific tool.
\par

\section{Evaluation of the perturbation solution} \label{s:averagedflow}

Once the perturbation solution have been computed by any means, it remains to be evaluated for given initial conditions. Recursions can obviously be efficiently programmed in platforms with limited memory. On the contrary, literal expression may require a much larger amount of memory to be stored. In view of memory limitation is not a major issue these days, one might be tempted to forgo the compact elegance of recursions in favor of the alacrity of explicit expressions. However, the use of recursions provides also a great versatility in the selection of the perturbation model, which is a non-negligible advantage in different applications. One important case in which this versatility is of great help is in the process of assessing the reliability of different simplified models to be used for mission design purposes.
\par

Thus, it commonly emerges the question on how many zonal harmonics of the gravitational potential must be kept into the perturbation model to capture the major features of the long-term dynamics while having a model as simple as possible to ease computations \cite{Lara2018Stardust}. The answer is normally ascertained after a process that includes the phase space representation for increasing complexity of the dynamical model. We will illustrate this procedure by exploring the mean elements phase space of high-inclination low-altitude lunar orbits for increasing number of zonal harmonics.
\par

The reduced flow stemming from Eq.~(\ref{Hmean}) is of 1-DOF and, for that reason, its orbits can be depicted by simple contour plots of the mean elements Hamiltonian without need of carrying out any numerical integration. Indeed, due to the axial symmetry of the zonal potential, Eq.~(\ref{zonalKaula}), the right ascension of the ascending node $\Omega=h$ is a cyclic variable. Therefore, its conjugate momentum $H$, the polar component of the angular momentum vector, is a dynamical constant. Besides, because the mean anomaly has been removed in the averaging procedure leading to Eq.~(\ref{Hmean}), its conjugate momentum $L$, the Delaunay action, is a formal dynamical constant of the mean elements Hamiltonian. In consequence, Eq.~(\ref{Hmean}) only depends on the remaining canonical variables, viz.~$g$ and $G$, and for given initial conditions remains constant because the averaged Hamiltonian does not depend explicitly on time.
\par

Alternatively, since $e=\sqrt{1-G^2/L^2}$ and $L$ is constant, on average, the mean elements Hamiltonian can be viewed as depending on the eccentricity and the argument of the periapsis $\omega=g$, in which the dynamical constants $a=L^2/\mu$ and $\sigma=H/L$ define a parameters plane. Furthermore, because $G=L$ for a circular orbit, $\sigma$ represents the cosine of the inclination of a circular orbit, and, hence, it is common to replace the dynamical parameter represented by $\sigma=\cos{I}_\mathrm{circular}$ by the corresponding inclination of circular orbits $I_\mathrm{circular}$ \cite{Lara2010,LaraSanJuanLopezOchoa2013a}. The visualization of the long-term dynamics by means of eccentricity vector diagrams, which are polar plots in which the eccentricity plays the role of the radius and the argument of the periapsis is the polar angle, is a useful tool for mission designers \cite{Cook1966,CuttingFrautnickBorn1978,Shapiro1996,FoltaQuinn2006}.
\par

\begin{figure}[htbp]\vspace{-1cm}
\centering
\includegraphics[scale=0.6]{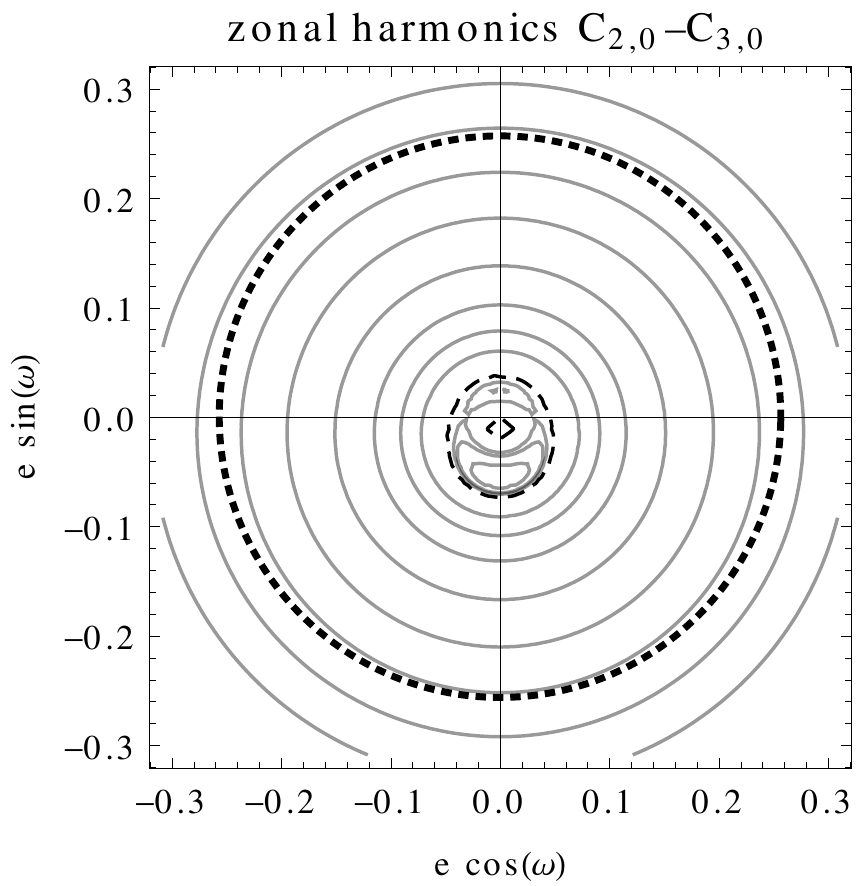} 
\includegraphics[scale=0.6]{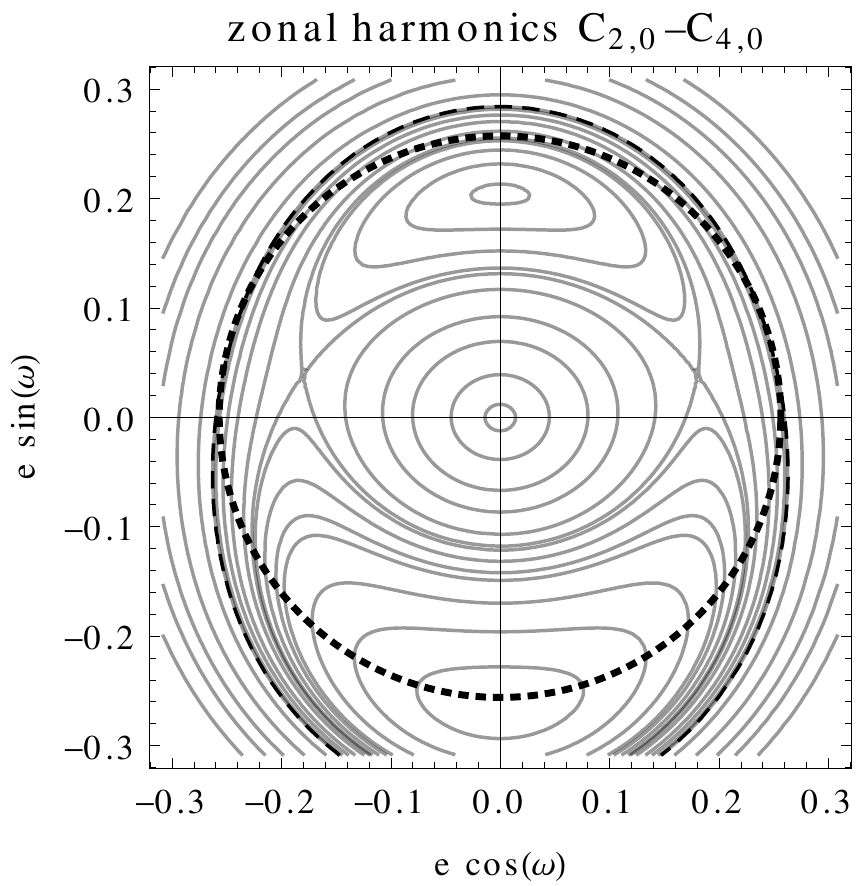} 
\includegraphics[scale=0.6]{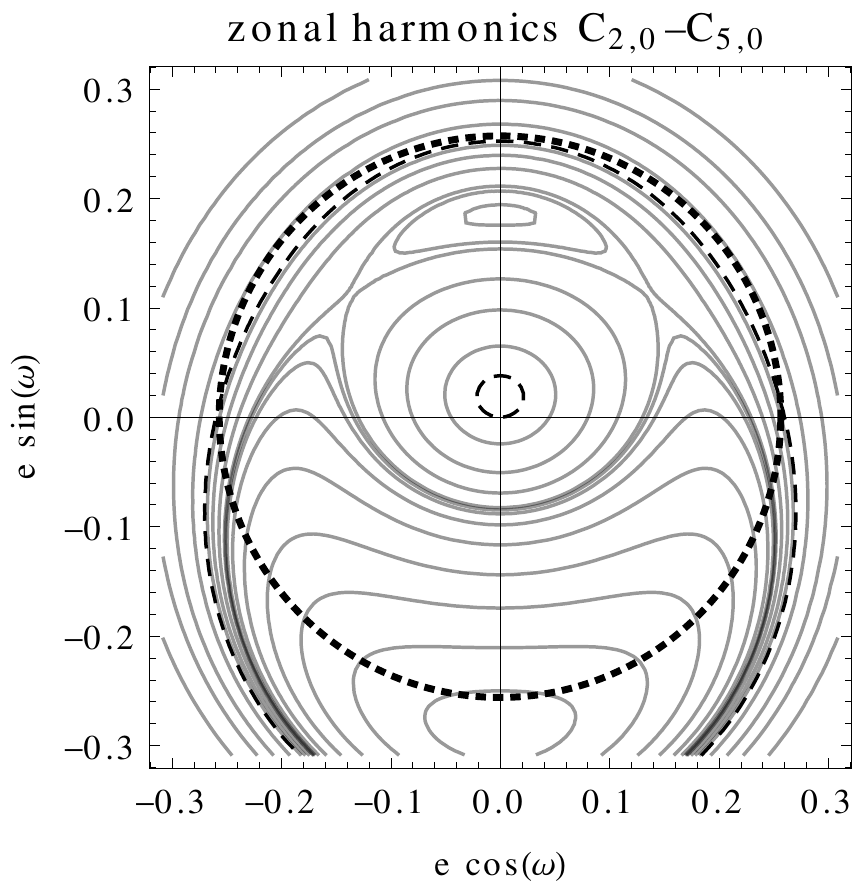} \\
\includegraphics[scale=0.6]{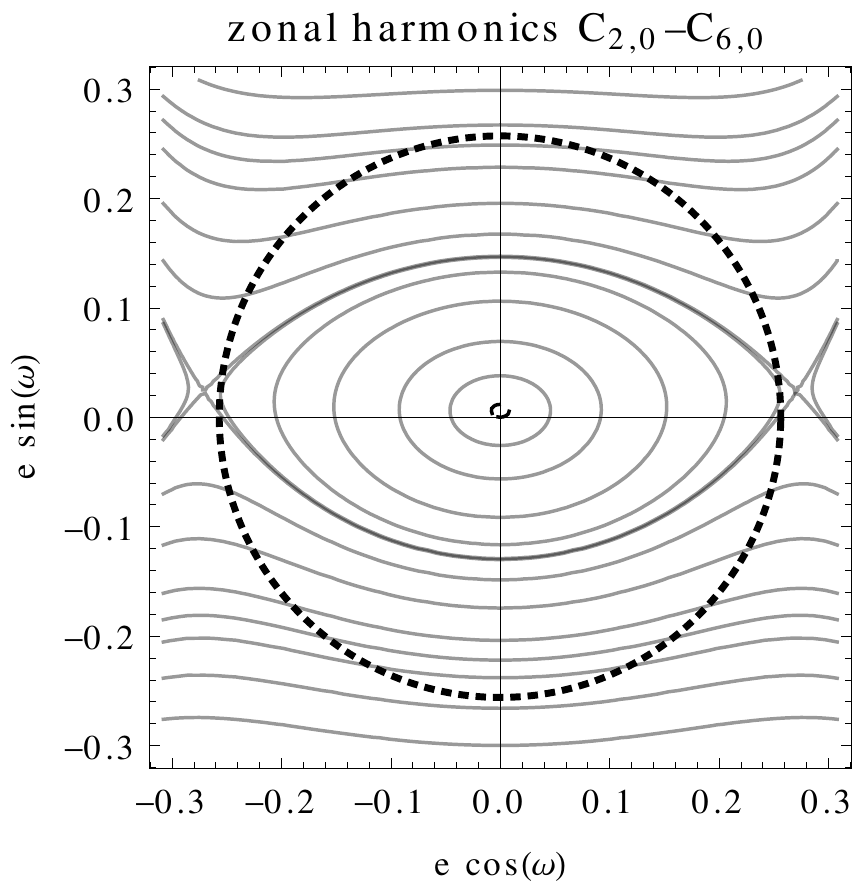}
\includegraphics[scale=0.6]{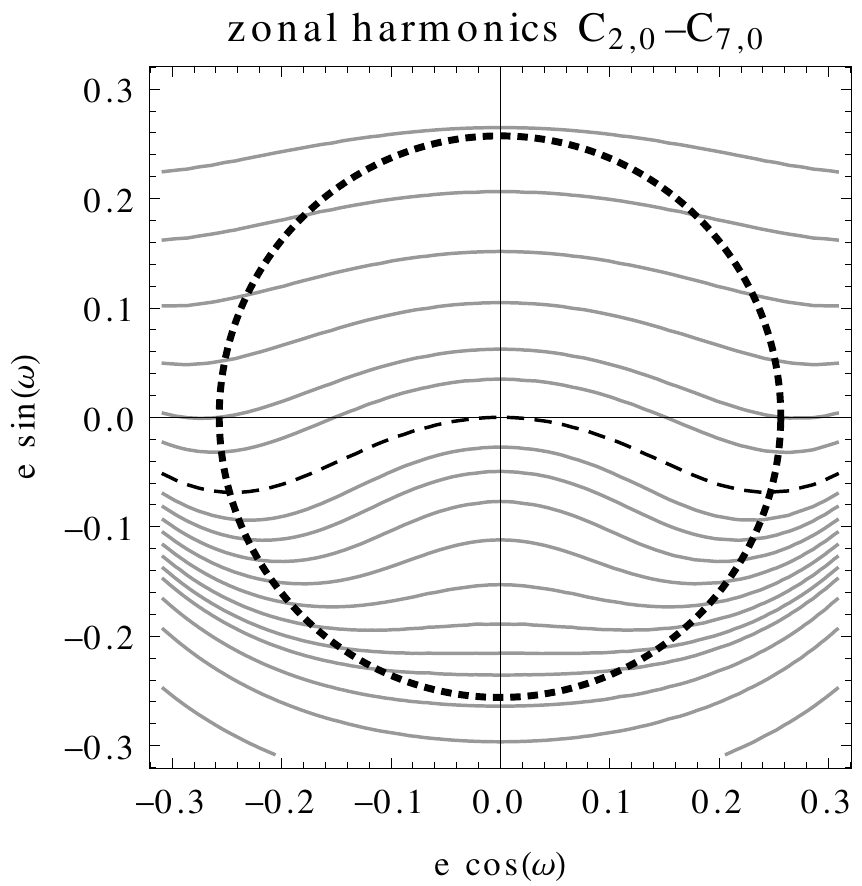}
\includegraphics[scale=0.6]{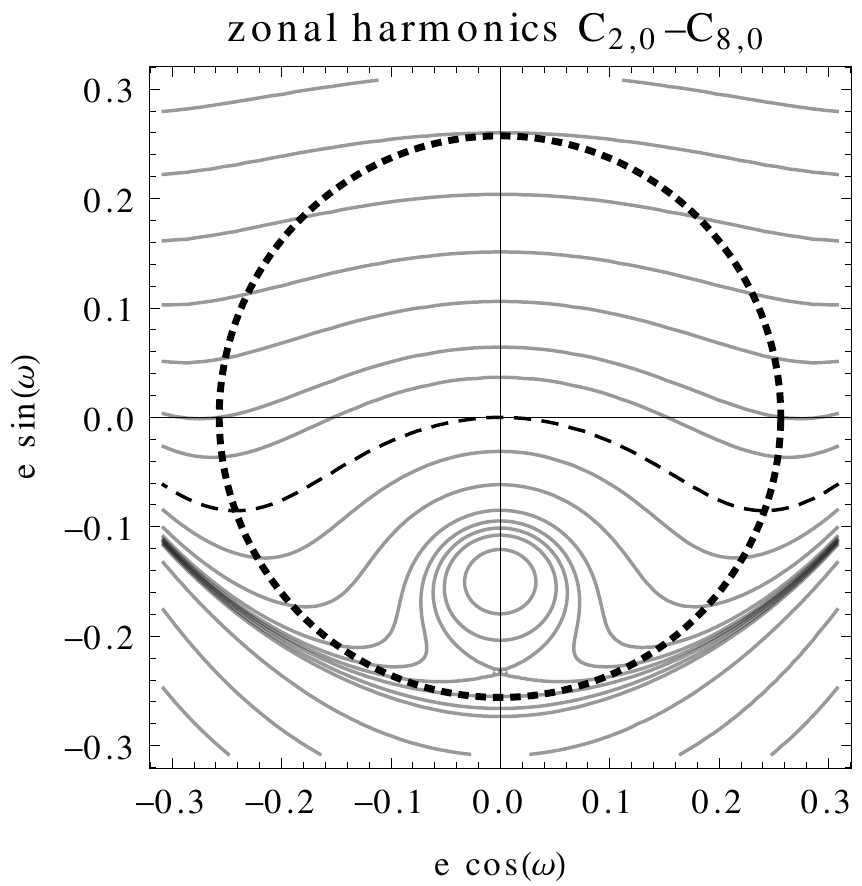} \\
\includegraphics[scale=0.6]{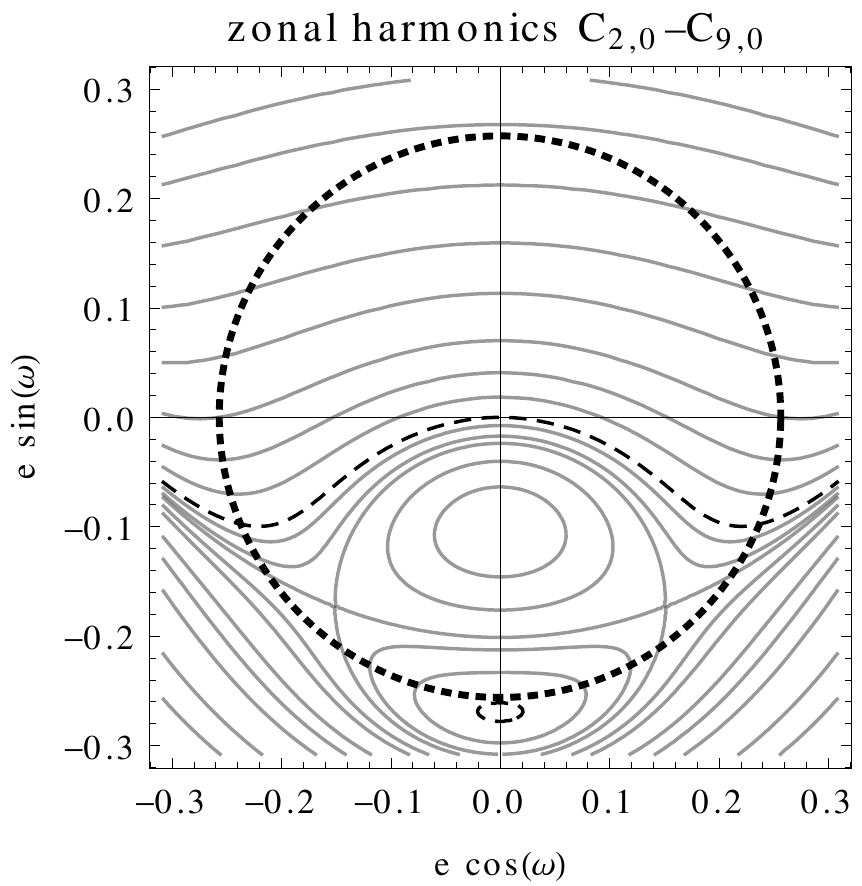}
\includegraphics[scale=0.6]{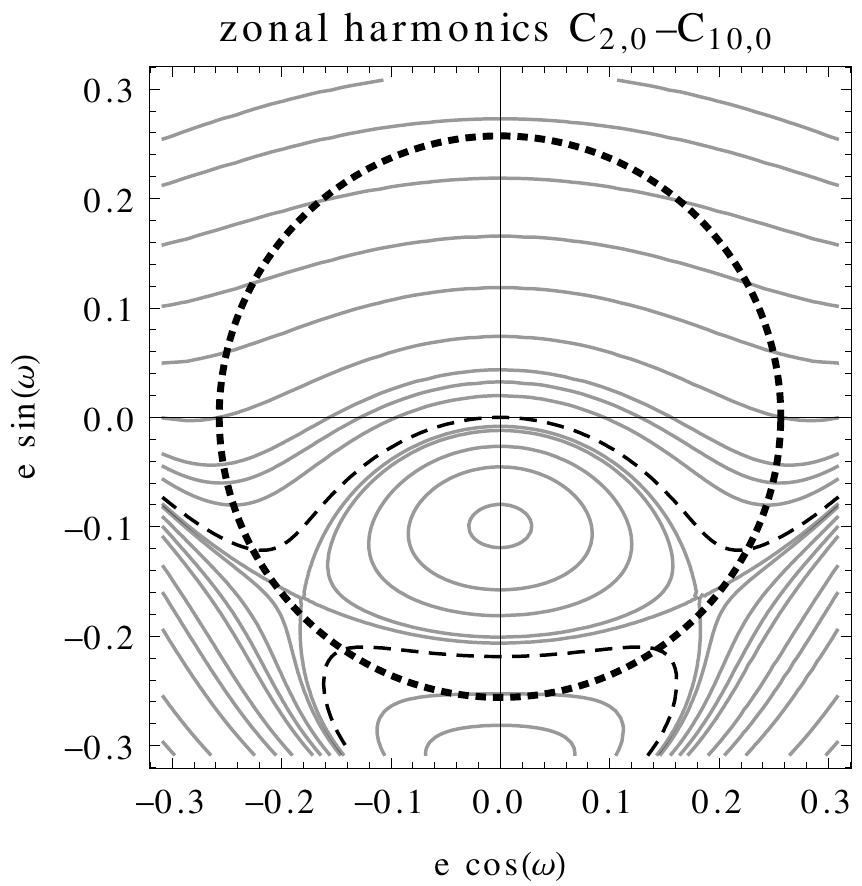}
\includegraphics[scale=0.6]{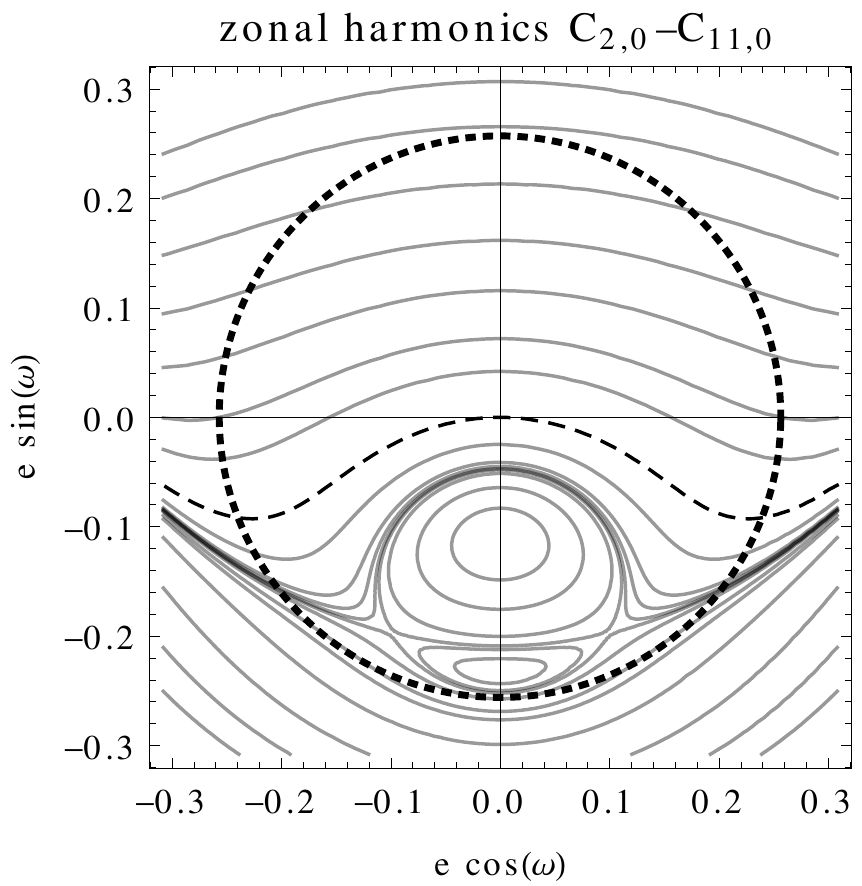} \\
\includegraphics[scale=0.6]{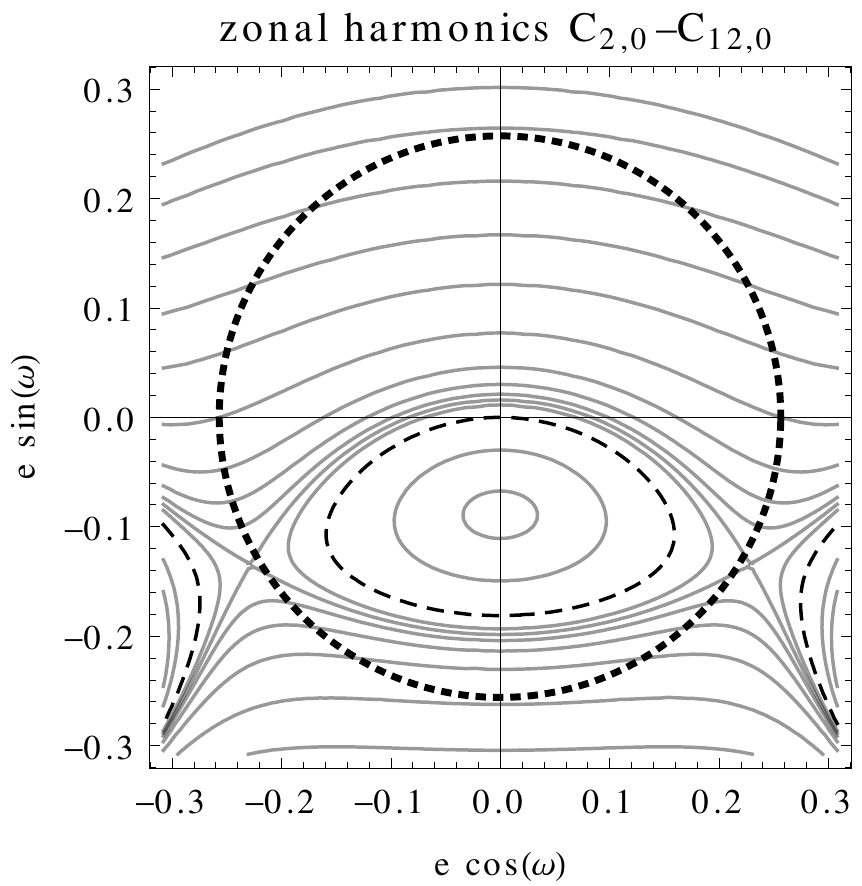}
\includegraphics[scale=0.6]{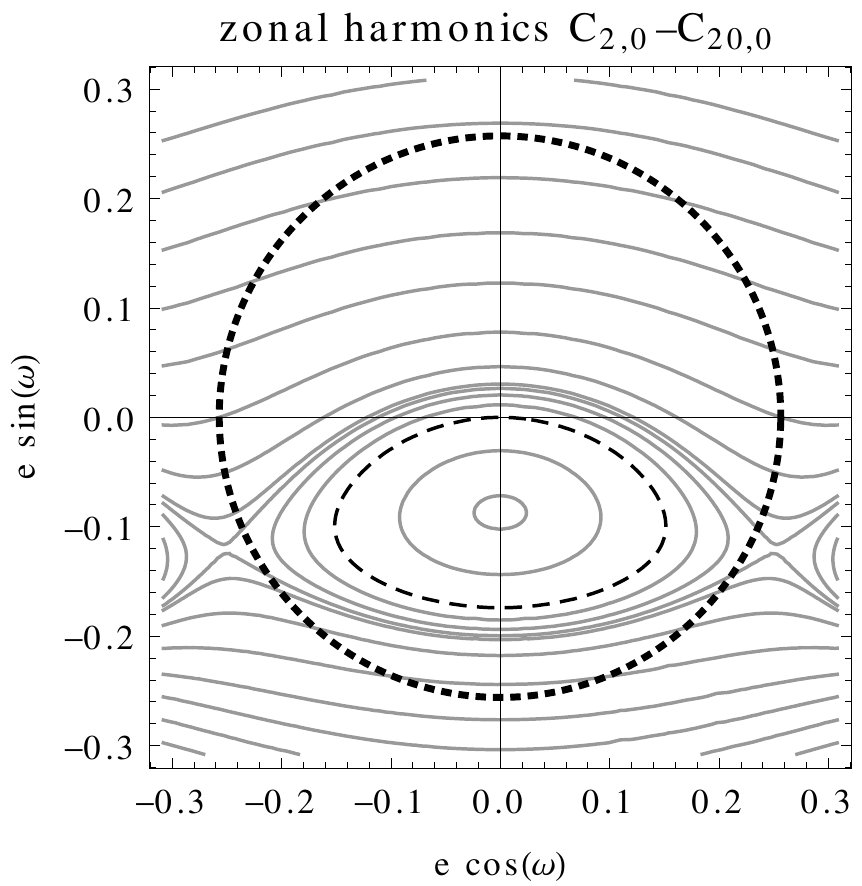}
\includegraphics[scale=0.6]{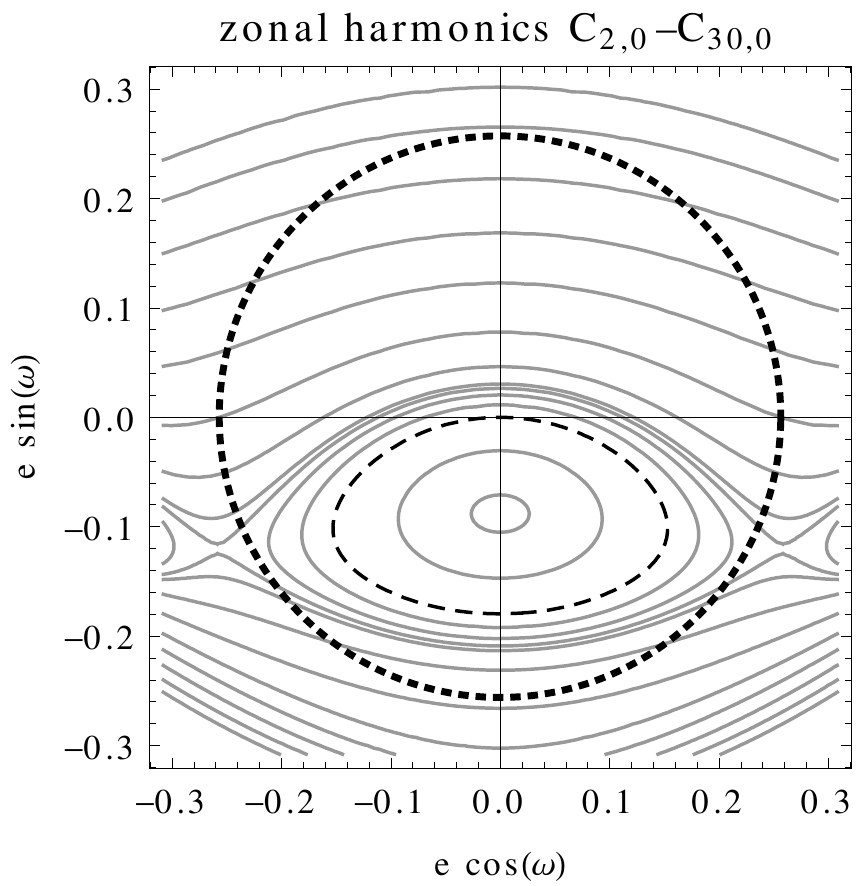}
\caption{Influence of the number of zonal harmonics in the long-term dynamics of a lunar orbit with $a=R_\oplus+600$ km over the surface of the moon, $I_\mathrm{circular}=63.45\deg$.
}
\label{f:ew600a6345i}
\end{figure}

\begin{figure}[htbp]
\centering
\includegraphics[scale=0.6]{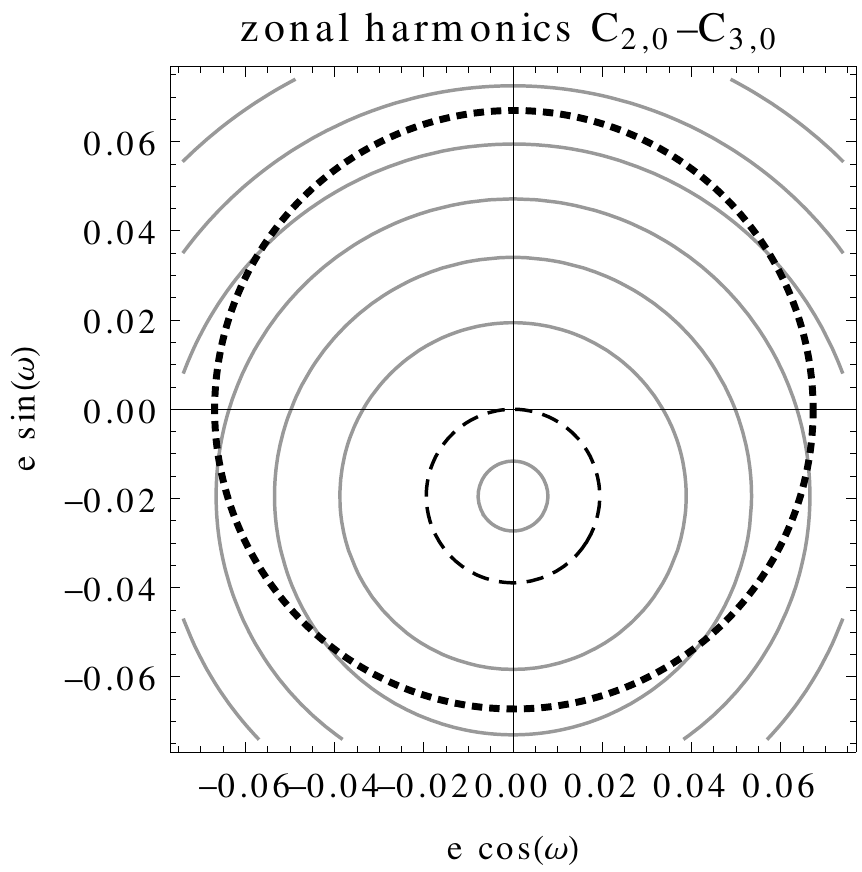} 
\includegraphics[scale=0.6]{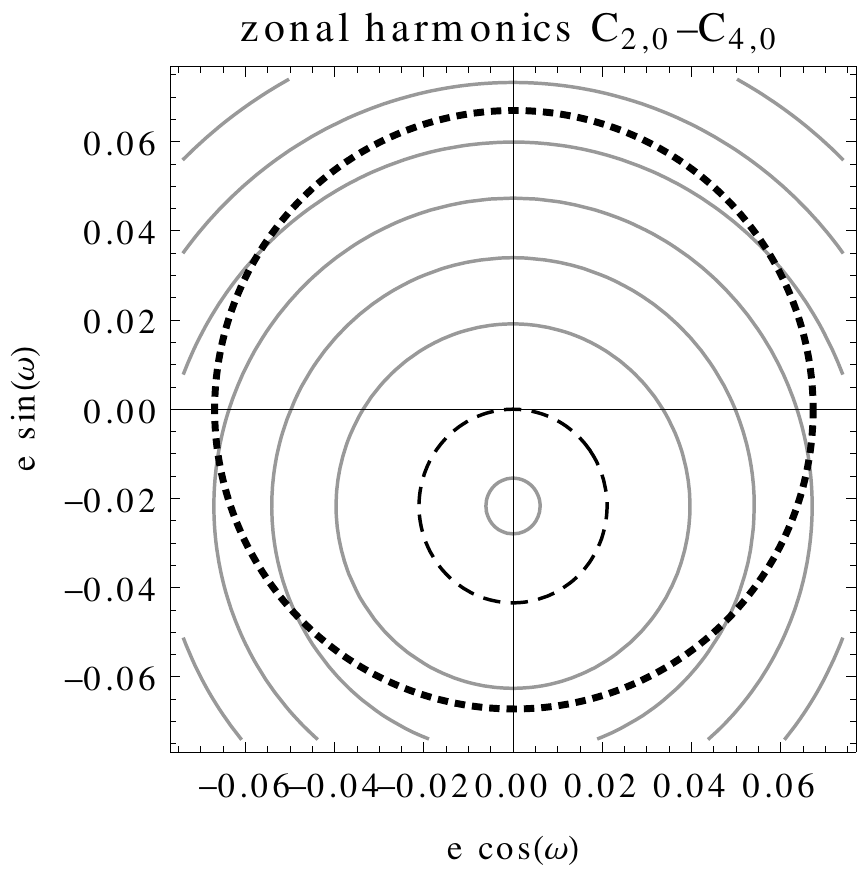} 
\includegraphics[scale=0.6]{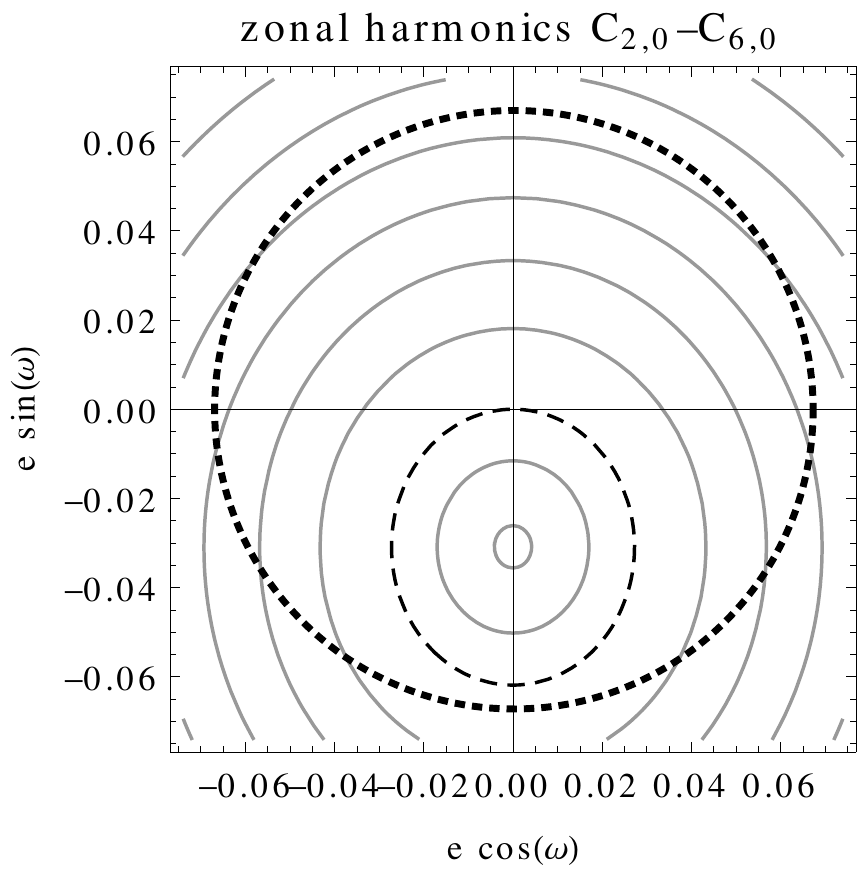} \\
\includegraphics[scale=0.6]{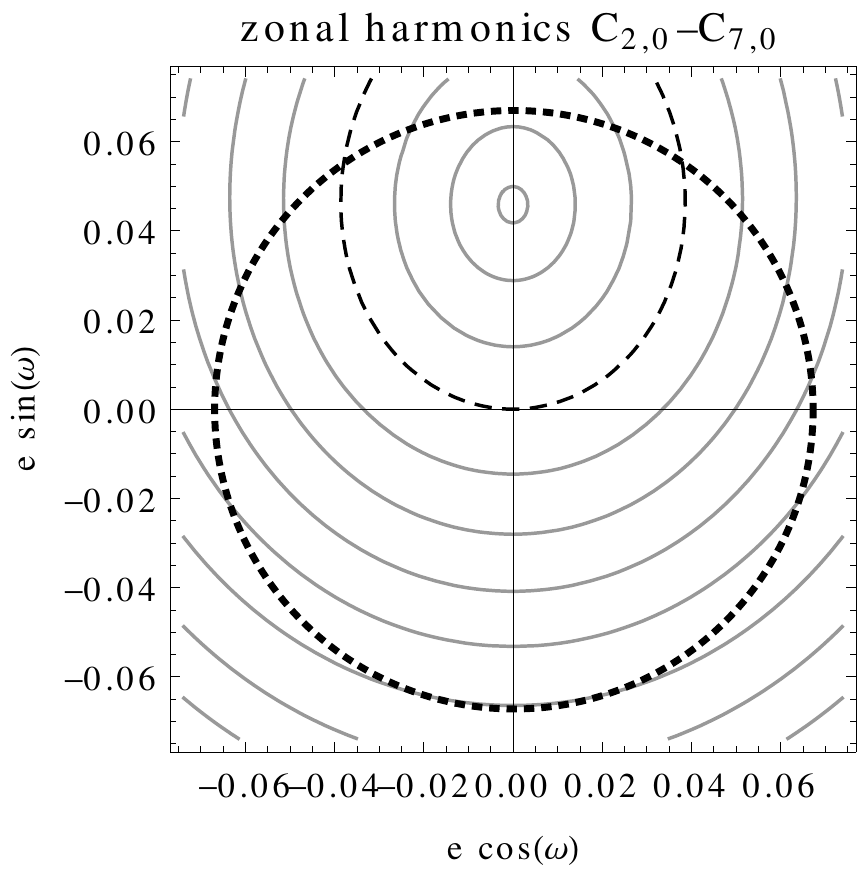}
\includegraphics[scale=0.6]{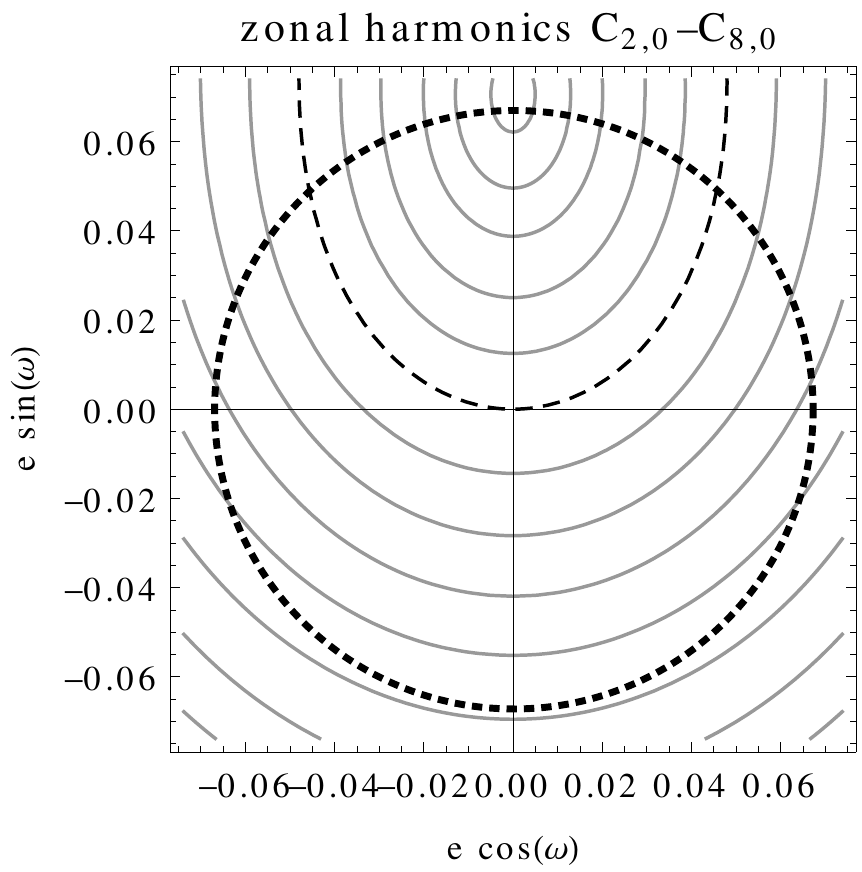}
\includegraphics[scale=0.6]{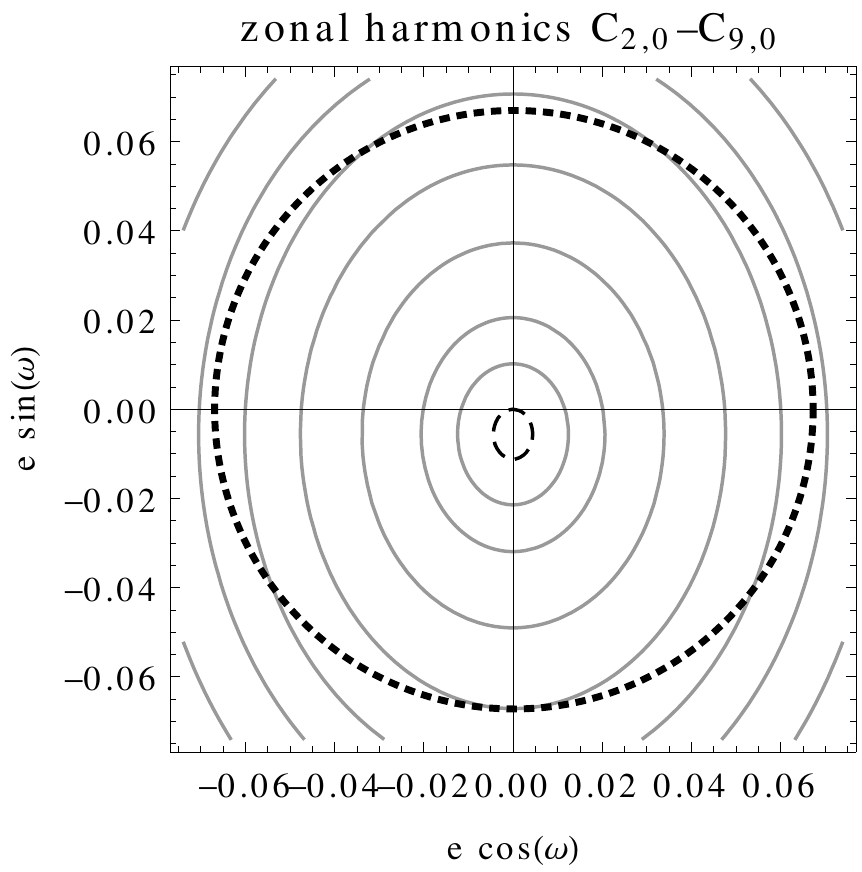} \\
\includegraphics[scale=0.6]{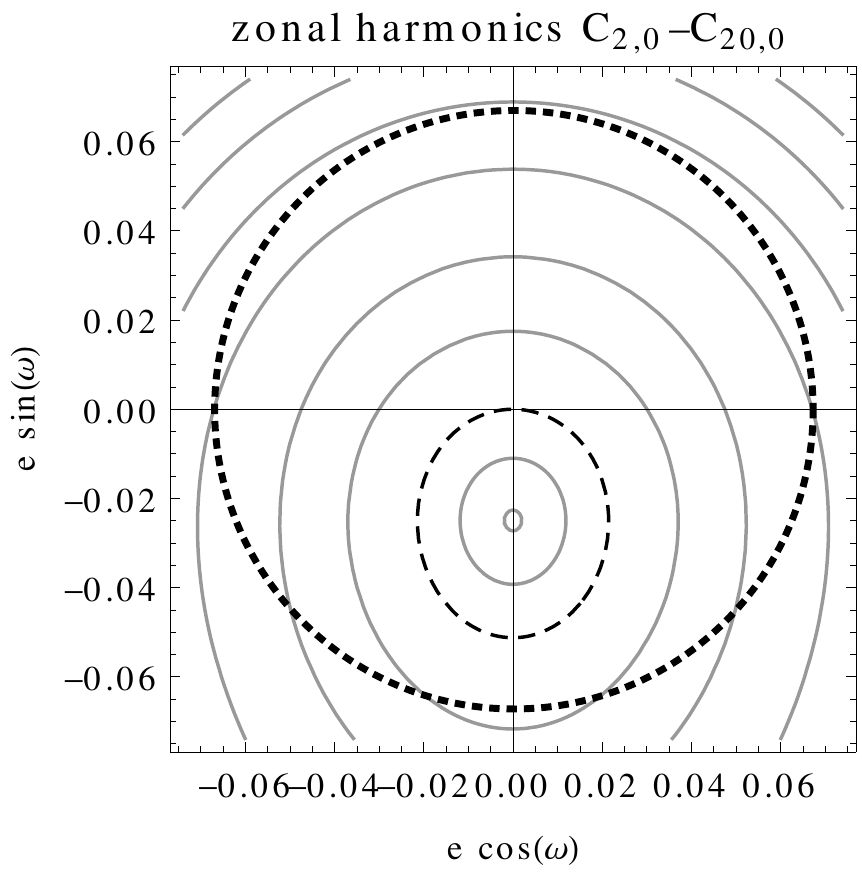}
\includegraphics[scale=0.6]{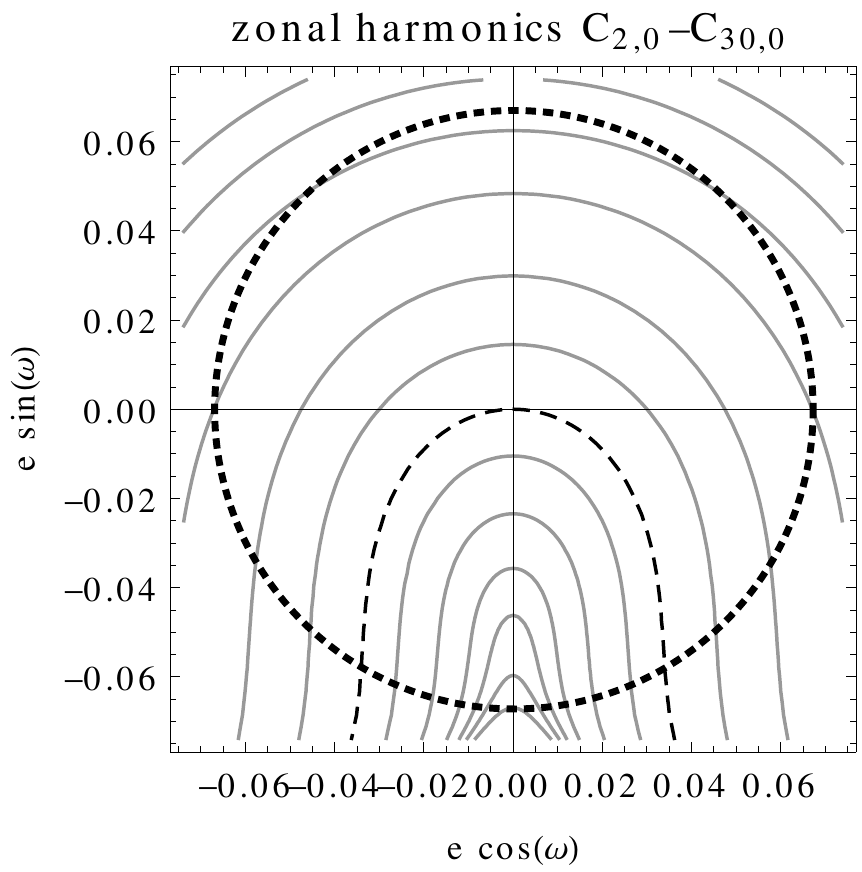}
\includegraphics[scale=0.6]{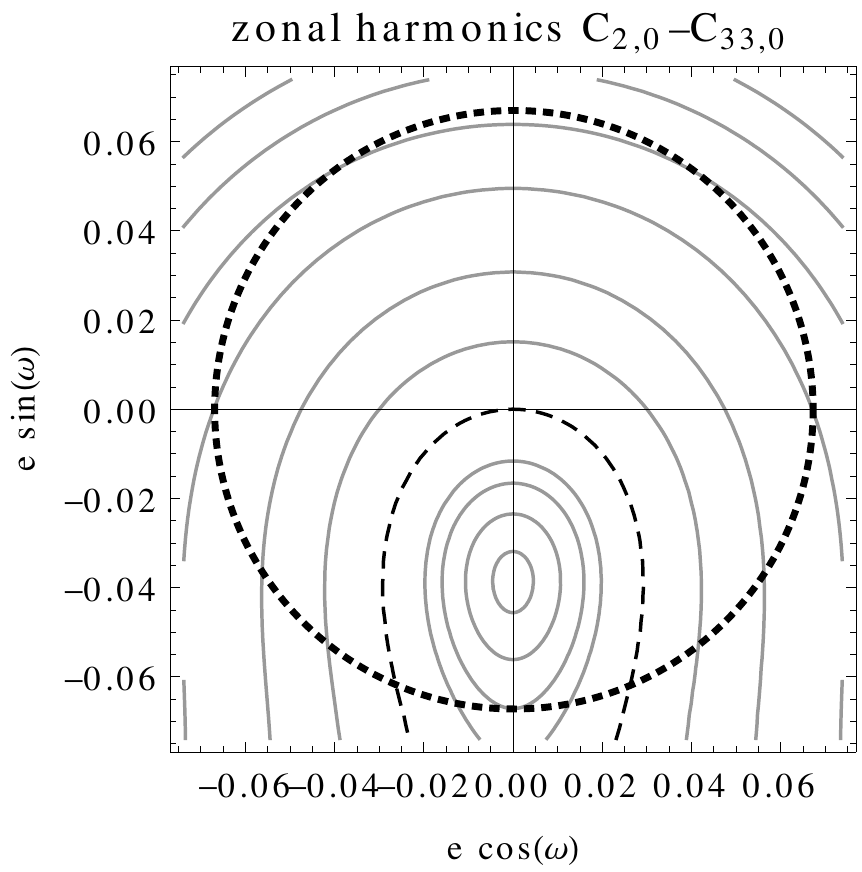}
\caption{Long-term dynamics of a lunar orbit with $a=R_\oplus+125$ km over the surface of the moon, $I_\mathrm{circular}=88\deg$.
}
\label{f:ew125km88deg}
\end{figure}

\begin{figure}[htb]
\centering
\includegraphics[scale=0.6]{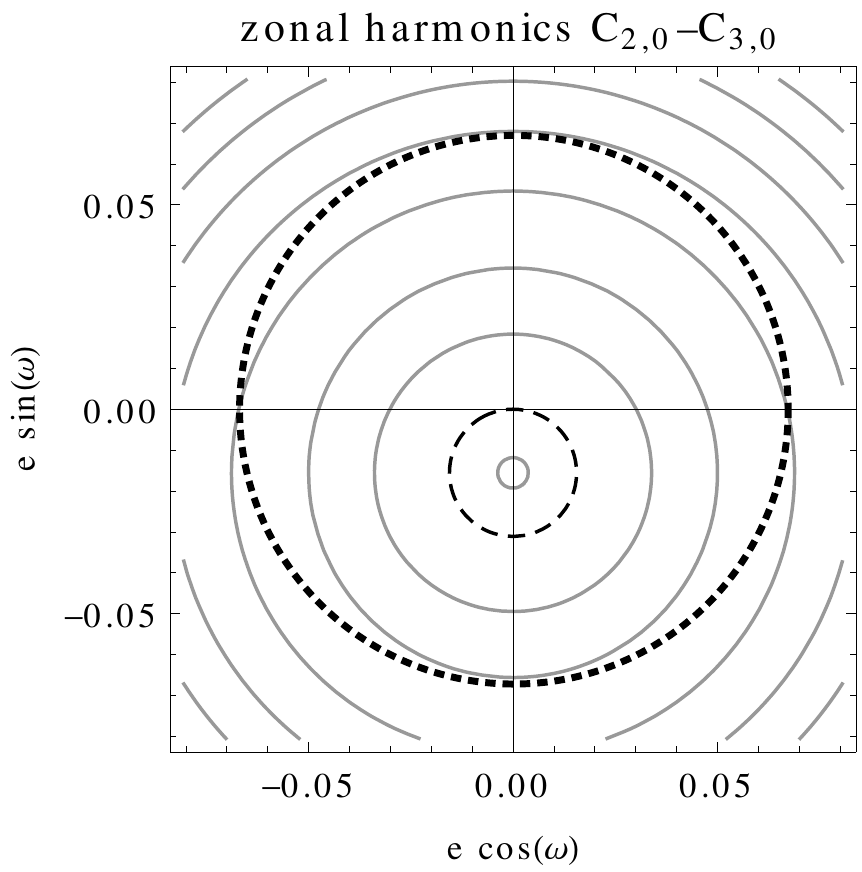} 
\includegraphics[scale=0.6]{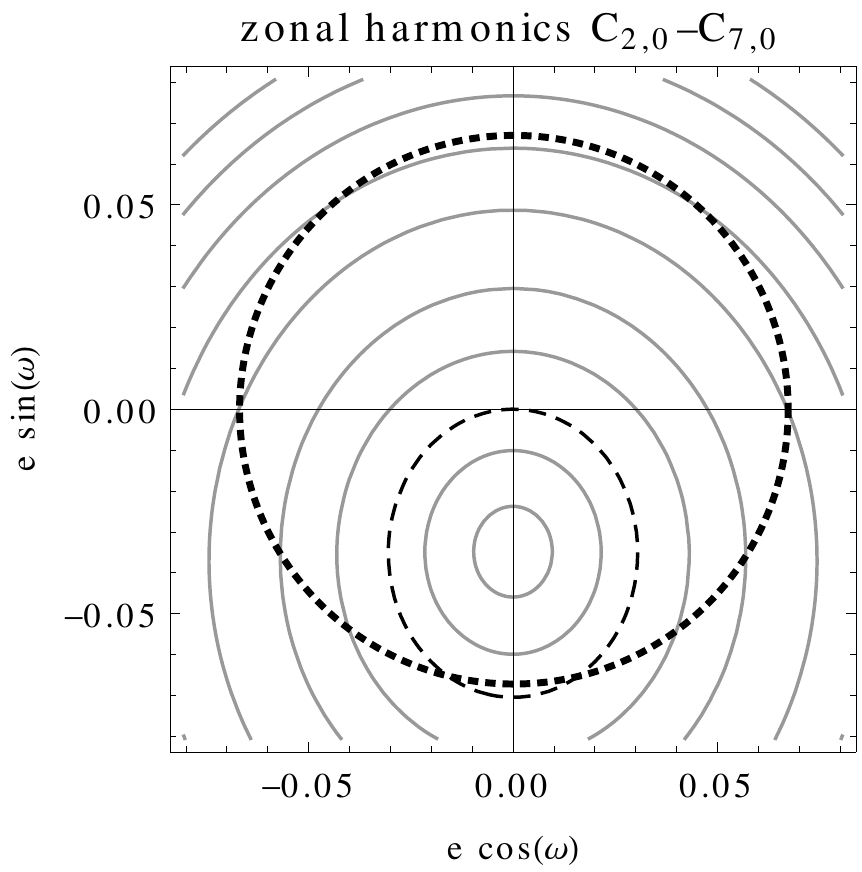} 
\includegraphics[scale=0.6]{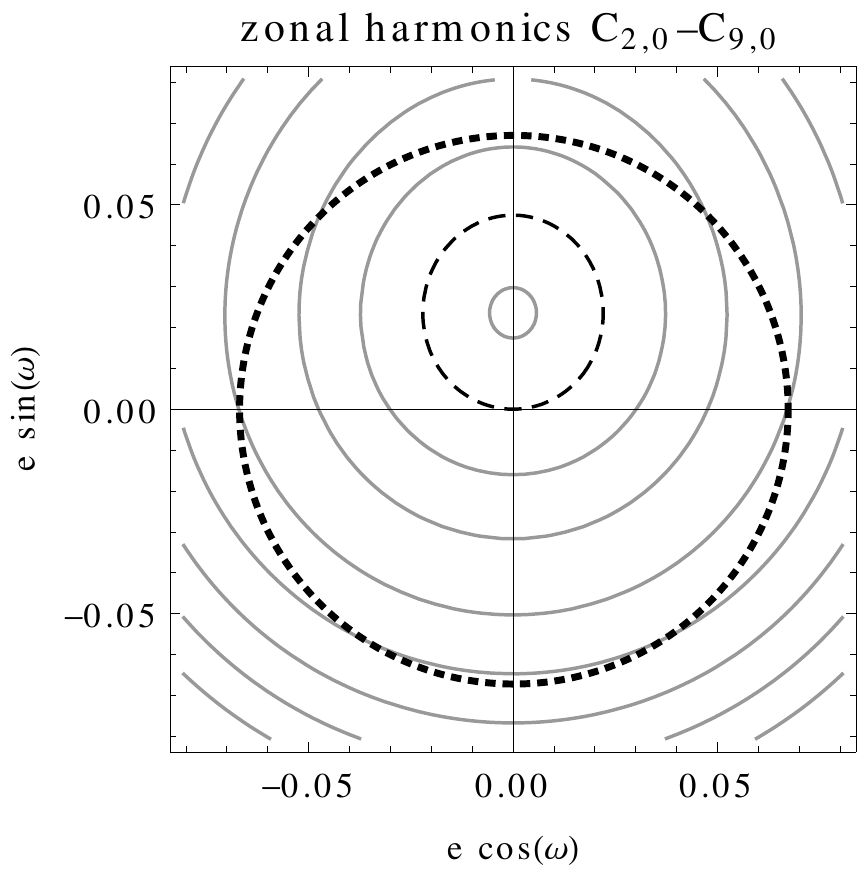} \\
\includegraphics[scale=0.6]{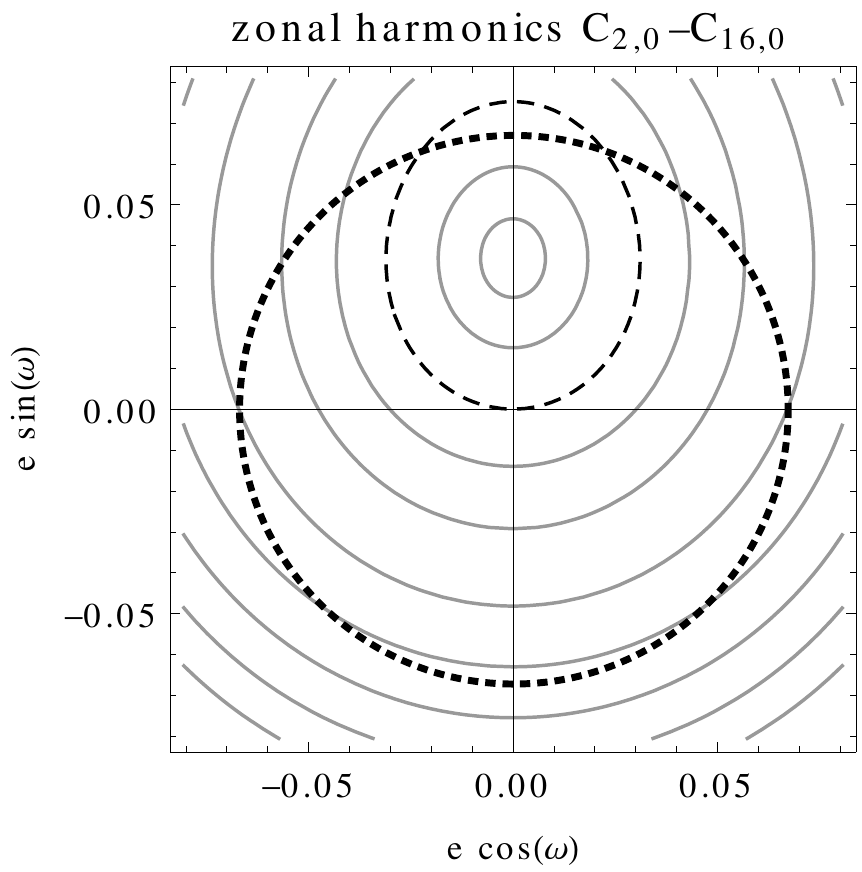}
\includegraphics[scale=0.6]{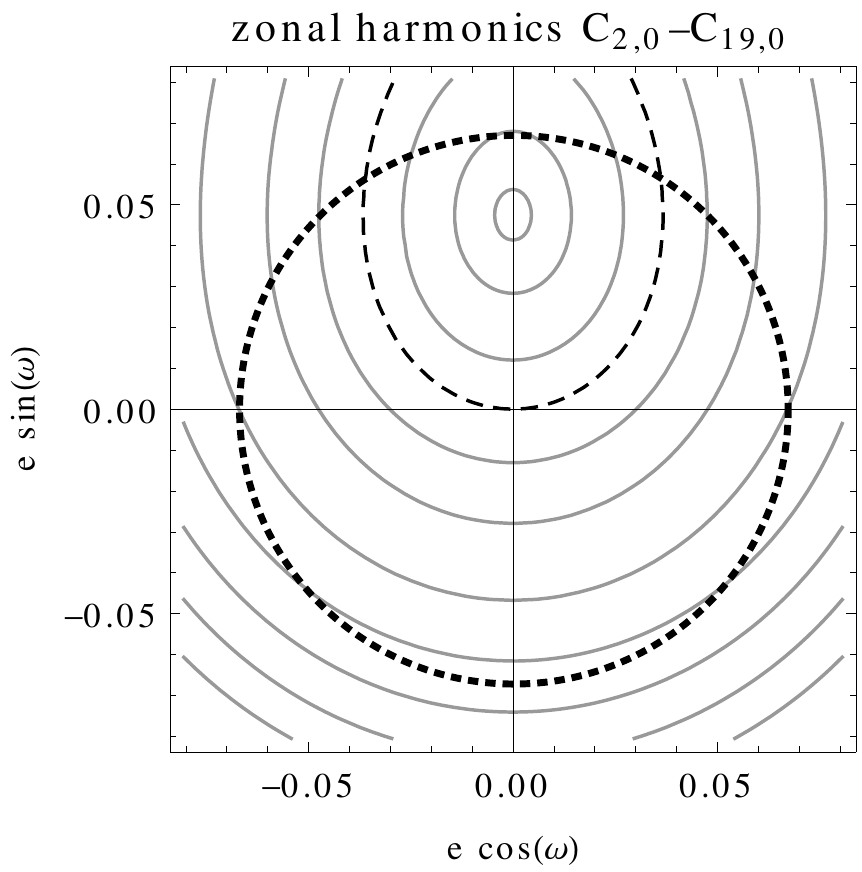}
\includegraphics[scale=0.6]{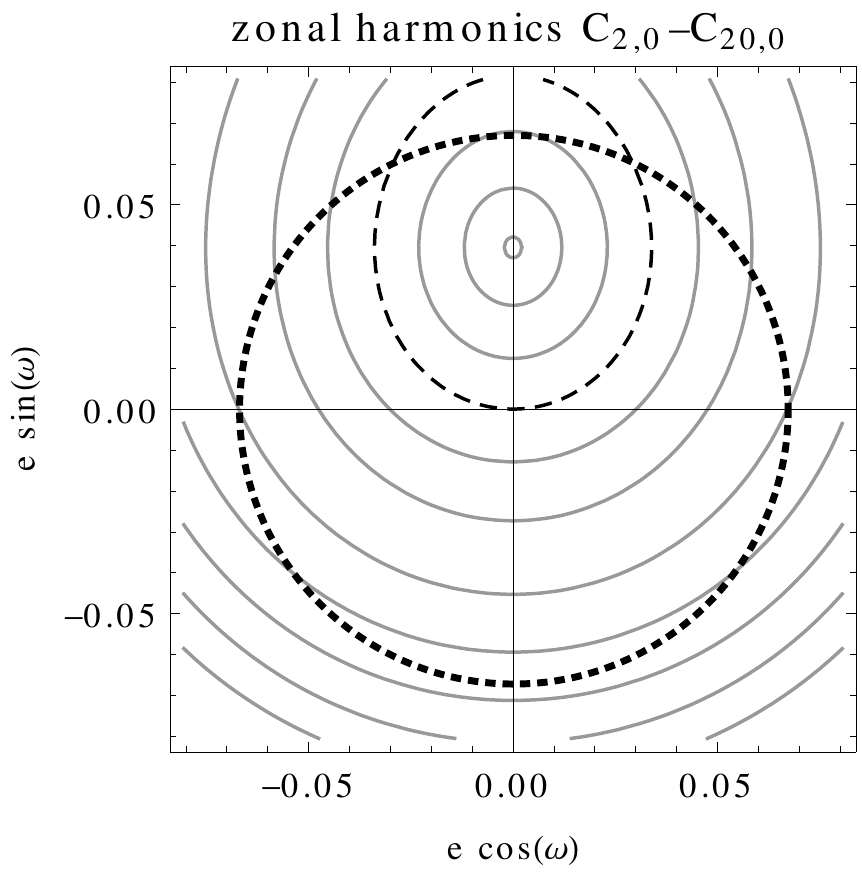}
\caption{Long-term dynamics of a lunar orbit with $a=R_\oplus+125$ km over the surface of the moon, $I_\mathrm{circular}=53\deg$.
}
\label{f:ew125km53deg}
\end{figure}

The differences in evaluating different truncations of the mean elements Hamiltonian for an orbit with desired characteristics is illustrated in Figs.~\ref{f:ew600a6345i}--\ref{f:ew125km53deg} for a lunar orbiter, in which the physical parameters of the Selenopotential have been taken form the lunar \texttt{lp150q} model \cite{Konoplivetal2001}. The dotted circumferences in the different plots of theses figures mark the eccentricity at which the orbits would impact with the surface of the moon. 
\par

Thus, Fig.~\ref{f:ew600a6345i} shows the case of a low lunar orbit in which the semi-major axis is about 1.3 times the lunar radius and the inclination of the circular orbits is fixed to 63.45 deg. As shown in the figure, the eccentricity vector diagrams that visualize the long-term dynamics for this point of the parametric plane notably change depending on the number of harmonics kept in the truncated model. Indeed, when only the Moon's $C_{2,0}$ and $C_{3,0}$ are taken into account, the perigee of non-impact orbits commonly circulates except in the case of almost circular orbits, where three stable frozen orbits exist, two of them with the argument of the perigee at $\omega=-\pi/2$, whereas the perigee of the other points to the opposite direction. Besides, two unstable almost circular frozen orbits exist with the arguments of the perigee pointing approximately to 0 and $\pi$. The simple addition of the moon's $C_{4,0}$ magnifies the eccentricity of four of the five previous frozen orbits, and only one of the stable orbits remain almost circular. The argument of the perigee of the unstable solutions clearly departs from the $0$--$\pi$ direction, and the eccentric stable frozen orbit with the argument of the perigee $\omega=-\pi/2$ becomes an impact orbit. The inclusion of more zonal harmonics in the mean elements Hamiltonian further modifies the figure, and when the moon's $C_{2,0}$--$C_{6,0}$ are taken into account, only the almost circular orbit survives as a non-impact orbit. Quite notably, no frozen non-impact orbits exist in the $C_{2,0}$--$C_{7,0}$ model, in which corresponding plot a dashed curve passing through the origin illustrates the evolution of the eccentricity of the circular orbits. Addition of more harmonics to the truncated hamiltonian show that the evolution of the circular orbits avoids impact for a $C_{2,0}$--$C_{12,0}$ model, a case in which a frozen non-impact exists with $e\approx0.09$ and $\omega=-\pi/2$. Finally, the situations seems to stabilize when $C_{2,0}$--$C_{30,0}$ are taken into account, and the addition of more zonal harmonics to the long-term Hamiltonian does not introduce relevant quantitative modifications of the long-term dynamics.
\par

The case of a low-altitude, high-inclination orbit is shown in Fig.~\ref{f:ew125km88deg}. Now, the parameters plane is set to the values $a=1.07$ lunar radius and $I_\mathrm{circular}=88$ deg, which would fit to the case of typical mapping orbits. Only one non-impact frozen orbit exists in this case, but the effect of the different zonal harmonics may change the argument of the perigee from $-\pi/2$ to $\pi/2$, or even prevent existence of non-impact frozen orbits. In this particular case, due to the lower altitude, the model needs a higher degree truncation for the long-term behavior to be stabilized, which happens when $C_{2,0}$--$C_{33,0}$ are taken into account, a case in which a non-impact frozen orbit exists whit $\omega=-\pi/2$ and $e\approx0.04$.
\par

Finally, the example shown in Fig.~\ref{f:ew125km53deg} shows that the long-term dynamics of lower-inclination orbits can be approached with lower degree truncations of the mean elements Hamiltonian. Thus, for an inclinations of the circular orbits $I_\mathrm{circular}=53$ deg the $C_{2,0}$--$C_{9,0}$ truncation suffices to capture the main characteristics of the dynamics, yet important quantitative variations are noted for higher degree truncations until the stabilization happens for the $C_{2,0}$--$C_{20,0}$ model.


\section{Conclusions}

It is well known that Kaula's seminal recursions provide an efficient way for the direct construction of the gravitational potential in orbital elements. Similar recursions derived from them can be used in the construction of the long-term gravitational potential, a case in which short-period effects are removed by averaging. Simulations in this paper show that Kaula's approach clearly remains as the benchmark to which the performance of other recursions in the literature, as well as the brut force approach, must be compared in the construction of a high degree long-term potential. On the other hand, second order effects of the second zonal harmonic, which fall out the scope of Kaula's approach, are needed in the investigation of the long-term propagation of orbits about the Earth or Earth-like bodies. However, these effects can be explicitly incorporated to the long-term disturbing function from available expressions in the literature. 
\par

The use of the averaged potential in the construction of eccentricity vector diagrams showed as an efficient way of ascertaining the complexity of the dynamical model to be used in the description of a particular region of phase space. It discloses the simplest model to be used in a real description of the long-term dynamics due to the non-centralities of the gravitational potential, yet additional disturbing effects, as for instance third-body perturbations, may be needed in the real description of the dynamics.

\subsection*{Acknowlwedgements}
Support by the Spanish State Research Agency and the European Regional Development Fund under Project ESP2016-76585-R (MINECO/ AEI/ERDF, EU) is recognized. ML also recognizes support from the Project ESP2017-87271-P of the same funding agencies.

\appendix

\section{Deprit's perturbation approach} \label{s:LieTransforms}

The philosophy of a perturbation method is to find a transformation of variables $\mathcal{T}:(\Vec\xi,\varepsilon)\rightarrow\Vec\xi'$, where $\Vec\xi,\Vec\xi'\in\mathbb{R}^m$ are given by an expansion in power series of the small parameter $\varepsilon$, viz.
\[
\Vec\xi=\Vec\xi'+\sum_{j\ge1}\frac{\varepsilon^j}{j!}\Vec{\xi}_j(\Vec\xi'),
\]
such that, after truncating the resulting series to some degree $\varepsilon^n$, the differential equations of the flow in the new variables are simpler to integrate than in the original ones. A common case is to find a transformation that, after truncation, decouples the differential system into a reduced one, which due to the reduced dimension should be simpler to solve than the original one, and a number of quadratures that can be integrated once the solution of the reduced system is known.
\par

For Hamiltonian problems the dimension $m$ is even, and the variables dissociate into $m/2$ coordinates $\Vec{x}$ and their corresponding conjugate momenta $\Vec{X}$. In that case, the transformation, which is now written $\mathcal{T}:(\Vec{x},\Vec{X},\varepsilon)\rightarrow(\Vec{x}',\Vec{X}')$, is derived from a single scalar generating function,
\begin{equation} \label{genfun}
\mathcal{W}=\sum_{m\ge0}\frac{\varepsilon^m}{m!}W_{m+1},
\end{equation}
which is conveniently computed using Deprit's perturbation method by Lie transforms \cite{Deprit1969}. Thus, starting from a Hamiltonian in the form of Eq.~(\ref{zonalHam}), in which $\mathcal{H}_{m,0}\equiv\mathcal{H}_{m,0}(\Vec{x},\Vec{X})$, the method provides a new Hamiltonian 
\begin{equation} \label{newHam}
\mathcal{K}=\sum_{m\ge0}\frac{\varepsilon^m}{m!}\mathcal{H}_{0,m}(\Vec{x}',\Vec{X}'),
\end{equation}
as well as the generating function (\ref{genfun}) by computing successive solutions of the homological equation
\begin{equation} \label{homological}
\{\mathcal{H}_{0,0};W_m\}+\widetilde{\mathcal{H}}_{0,m}=\mathcal{H}_{0,m},
\end{equation}
in which
\begin{itemize}
\item terms $\widetilde{\mathcal{H}}_{0,m}$ are known from previous computations;
\item terms $\mathcal{H}_{0,m}$ are chosen at will, in accordance with the simplification criterion;
\item terms $W_m$ are solved from the partial differential equation resulting from the evaluation of the Poisson bracket $\{\mathcal{H}_{0,0};W_m\}$ ---sometimes called the Lie operator.
\end{itemize}

At each step $m$, the homological equation (\ref{homological}) is obtained after rearranging Deprit's famous recursion
\begin{equation} \label{deprittriangle}
\mathcal{H}_{n,q+1}=\mathcal{H}_{n+1,q}+\sum_{0\le{m}\le{n}}{n\choose{m}}\,\{\mathcal{H}_{n-m,q};W_{m+1}\}.
\end{equation}
where the computation of terms $\mathcal{H}_{0,n}$ require the previous computation of all the terms $\mathcal{H}_{i,n-i}$, $0<i<n$, from the same recursion. For instance: $n=0$, $q=1$, leads to the sequence
\begin{eqnarray*}
\mathcal{H}_{0,2} &=&\mathcal{H}_{1,1}+\{\mathcal{H}_{0,1};W_{1}\} \\
\mathcal{H}_{1,1} &=& 
\mathcal{H}_{2,0}+\{\mathcal{H}_{1,0};W_{1}\}+\{\mathcal{H}_{0,0};W_{2}\}, \\
\mathcal{H}_{0,1} &=& \mathcal{H}_{1,0}+\{\mathcal{H}_{0,0};W_{1}\}.
\end{eqnarray*}
Hence, from Eq.~(\ref{homological}),
\[
\widetilde{\mathcal{H}}_{0,2}=\{\mathcal{H}_{0,1},W_1\}+\{\mathcal{H}_{1,0},W_1\}+\mathcal{H}_{2,0}.
\]
\par

Recursion (\ref{deprittriangle}) is not constrained to the computation of the new Hamiltonian and applies to any function $F=\sum_{m\ge0}(\varepsilon^m/m!)F_{m,0}(\Vec{x},\Vec{X})$. Hence, once $\mathcal{W}$ has been solved up to the desired order, Eq.~(\ref{deprittriangle}) can be used to compute the transformation equations that, up to the truncation order, transform $\mathcal{H}$ into $\mathcal{K}$.

For perturbed Keplerian motion, the homological equation (\ref{homological}) is conveniently solved in Delaunay canonical variables. In that case $\{\mathcal{H}_{0,0};W_m\}=-n{\partial{W}_m}/{\partial{\ell}}$, and the homological equation
\begin{equation} \label{homoDelo}
-n\frac{\partial{W}_m}{\partial{\ell}}+\widetilde{\mathcal{H}}_{0,m}=\mathcal{H}_{0,m},
\end{equation}
can be solved by quadrature.


\begin{thebibliography}{10}
\expandafter\ifx\csname url\endcsname\relax
  \def\url#1{\texttt{#1}}\fi
\expandafter\ifx\csname urlprefix\endcsname\relax\def\urlprefix{URL }\fi
\expandafter\ifx\csname href\endcsname\relax
  \def\href#1#2{#2} \def\path#1{#1}\fi

\bibitem{Kozai1963}
Y.~Kozai, Motion of a lunar orbiter, Publications of the Astronomical Society
  of Japan 15 (1963) 301--312.

\bibitem{CuttingFrautnickBorn1978}
E.~{Cutting}, J.~C. {Frautnick}, G.~H. {Born}, {Orbit analysis for Seasat-A},
  Journal of the Astronautical Sciences 26 (1978) 315--342.

\bibitem{ScheeresGumanVillac2001}
D.~J. {Scheeres}, M.~D. {Guman}, B.~F. {Villac}, {Stability Analysis of
  Planetary Satellite Orbiters: Application to the Europa Orbiter}, Journal of
  Guidance Control Dynamics 24~(4) (2001) 778--787.
\newblock \href {http://dx.doi.org/10.2514/2.4778} {\path{doi:10.2514/2.4778}}.

\bibitem{LaraSanJuan2005}
M.~Lara, J.~San-Juan, {Dynamic Behavior of an Orbiter Around Europa}, Journal
  of Guidance, Control and Dynamics 28~(2) (2005) 291--297.
\newblock \href {http://dx.doi.org/10.2514/1.5686} {\path{doi:10.2514/1.5686}}.

\bibitem{ColomboLuckingMcInnes2012}
C.~{Colombo}, C.~{L\"ucking}, C.~R. {McInnes},
  \href{http://www.sciencedirect.com/science/article/pii/S009457651200272X}{{Orbital
  dynamics of high area-to-mass ratio spacecraft with $J_2$ and solar radiation
  pressure for novel Earth observation and communication services}}, Acta
  Astronautica 81~(1) (2012) 137 -- 150.
\newblock \href
  {http://dx.doi.org/https://doi.org/10.1016/j.actaastro.2012.07.009}
  {\path{doi:https://doi.org/10.1016/j.actaastro.2012.07.009}}.
\newline\urlprefix\url{http://www.sciencedirect.com/science/article/pii/S009457651200272X}

\bibitem{ArmellinSanJuanLara2015}
R.~Armellin, J.~F. San-Juan, M.~Lara,
  \href{http://www.sciencedirect.com/science/article/pii/S0273117715002227}{End-of-life
  disposal of high elliptical orbit missions: The case of {INTEGRAL}}, Advances
  in Space Research 56~(3) (2015) 479--493, {Advances in Asteroid and Space
  Debris Science and Technology - Part 1}.
\newblock \href {http://dx.doi.org/10.1016/j.asr.2015.03.020}
  {\path{doi:10.1016/j.asr.2015.03.020}}.
\newline\urlprefix\url{http://www.sciencedirect.com/science/article/pii/S0273117715002227}

\bibitem{Colomboetal2015}
C.~Colombo, E.~M. Alessi, W.~van~der Weg, S.~Soldini, F.~Letizia, M.~Vetrisano,
  M.~Vasile, A.~Rossi, M.~Landgraf,
  \href{http://www.sciencedirect.com/science/article/pii/S0094576514004342}{End-of-life
  disposal concepts for libration point orbit and highly elliptical orbit
  missions}, Acta Astronautica 110 (2015) 298 -- 312, dynamics and Control of
  Space Systems.
\newblock \href
  {http://dx.doi.org/http://dx.doi.org/10.1016/j.actaastro.2014.11.002}
  {\path{doi:http://dx.doi.org/10.1016/j.actaastro.2014.11.002}}.
\newline\urlprefix\url{http://www.sciencedirect.com/science/article/pii/S0094576514004342}

\bibitem{Broucke1994}
R.~A. {Broucke}, {Numerical integration of periodic orbits in the main problem
  of artificial satellite theory}, Celestial Mechanics and Dynamical Astronomy
  58~(2) (1994) 99--123.
\newblock \href {http://dx.doi.org/10.1007/BF00695787}
  {\path{doi:10.1007/BF00695787}}.

\bibitem{JorbaMasdemont1999}
{\`A}.~{Jorba}, J.~{Masdemont}, {Dynamics in the center manifold of the
  collinear points of the restricted three body problem}, Physica D Nonlinear
  Phenomena 132 (1999) 189--213.
\newblock \href {http://dx.doi.org/10.1016/S0167-2789(99)00042-1}
  {\path{doi:10.1016/S0167-2789(99)00042-1}}.

\bibitem{Lara1999}
M.~Lara, {Searching for Repeating Ground Track Orbits: A Systematic Approach},
  Journal of the Astronautical Sciences 47~(3-4) (1999) 177--188.

\bibitem{GomezMondelo2001}
G.~G\'omez, J.~Mondelo,
  \href{http://www.sciencedirect.com/science/article/pii/S0167278901003128}{The
  dynamics around the collinear equilibrium points of the {RTBP}}, Physica D:
  Nonlinear Phenomena 157~(4) (2001) 283 -- 321.
\newblock \href
  {http://dx.doi.org/https://doi.org/10.1016/S0167-2789(01)00312-8}
  {\path{doi:https://doi.org/10.1016/S0167-2789(01)00312-8}}.
\newline\urlprefix\url{http://www.sciencedirect.com/science/article/pii/S0167278901003128}

\bibitem{Lara2003}
M.~{Lara}, {Repeat Ground Track Orbits of the Earth Tesseral Problem as
  Bifurcations of the Equatorial Family of Periodic Orbits}, Celestial
  Mechanics and Dynamical Astronomy 86~(2) (2003) 143--162.

\bibitem{LaraRussellVillac2007}
M.~{Lara}, R.~{Russell}, B.~F. {Villac}, {Classification of the Distant
  Stability Regions at Europa}, Journal of Guidance Control Dynamics 30 (2007)
  409--418.
\newblock \href {http://dx.doi.org/10.2514/1.22372}
  {\path{doi:10.2514/1.22372}}.

\bibitem{RussellLara2007}
R.~Russell, M.~Lara, Long-lifetime lunar repeat ground track orbits, Journal of
  Guidance, Control, and Dynamics 30~(4) (2007) 982--993.

\bibitem{LaraRussellVillac2006}
M.~{Lara}, R.~P. {Russell}, B.~{Villac},
  \href{https://ntrs.nasa.gov/search.jsp?R=20060042732}{{Stability Maps, Global
  Dynamics and Transfers (AAS 05-378)}}, in: B.~G. Williams, L.~A. D'Amario,
  K.~C. Howell, F.~R. Hoots (Eds.), AAS/AIAA Astrodynamics 2005, Vol. 123 of
  Advances in the Astronautical Sciences, American Astronautical Society,
  Univelt, Inc., P.O.~Box 28130, San Diego, California 92198, USA, 2006, pp.
  1983--2002.
\newline\urlprefix\url{https://ntrs.nasa.gov/search.jsp?R=20060042732}

\bibitem{LaraRussellVillac2007Meccanica}
M.~{Lara}, R.~{Russell}, B.~F. {Villac},
  \href{https://link.springer.com/article/10.1007/s11012-007-9060-z}{{Fast
  estimation of stable regions in real models}}, Meccanica 42~(5) (2007)
  511--515.
\newblock \href {http://dx.doi.org/10.1007/s11012-007-9060-z}
  {\path{doi:10.1007/s11012-007-9060-z}}.
\newline\urlprefix\url{https://link.springer.com/article/10.1007/s11012-007-9060-z}

\bibitem{Villac2008}
B.~F. {Villac}, {Using FLI maps for preliminary spacecraft trajectory design in
  multi-body environments}, Celestial Mechanics and Dynamical Astronomy 102
  (2008) 29--48.
\newblock \href {http://dx.doi.org/10.1007/s10569-008-9158-1}
  {\path{doi:10.1007/s10569-008-9158-1}}.

\bibitem{CoffeyDepritDeprit1994}
S.~L. {Coffey}, A.~{Deprit}, E.~{Deprit}, {Frozen orbits for satellites close
  to an earth-like planet}, Celestial Mechanics and Dynamical Astronomy 59~(1)
  (1994) 37--72.
\newblock \href {http://dx.doi.org/10.1007/BF00691970}
  {\path{doi:10.1007/BF00691970}}.

\bibitem{SanJuanLaraFerrer2006}
J.~F. {San-Juan}, M.~{Lara}, S.~{Ferrer}, {Phase Space Structure Around Oblate
  Planetary Satellites}, Journal of Guidance Control Dynamics 29 (2006)
  113--120.
\newblock \href {http://dx.doi.org/10.2514/1.13385}
  {\path{doi:10.2514/1.13385}}.

\bibitem{Konoplivetal1994}
A.~Konopliv, W.~Sjogren, R.~Wimberly, R.~Cook, A.~Vijayaraghavan, {A High
  Resolution Lunar Gravity Field and Predicted Orbit Behavior (AAS 93-622)},
  in: Advances in the astronautical Sciences, Vol.~85, American Astronautical
  Society, Univelt, Inc., P.O.~Box 28130, San Diego, California 92198, USA,
  1994, pp. 1275--1295.

\bibitem{Roncoli2005}
R.~Roncoli, {Lunar Constants and Models Document}, Tech. Rep. JPL D-32296, Jet
  Propulsion Laboratory, California Institute of Technology (Sep. 2005).

\bibitem{LaraSaedeleerFerrer2009}
M.~Lara, B.~de~Saedeleer, S.~Ferrer,
  \href{http://issfd.org/ISSFD_2009/InterMissionDesignII/Lara.pdf}{{Preliminary
  design of low lunar orbits}}, in: Proceedings of the 21st International
  Symposium on Space Flight Dynamics, Toulouse, France, ISSFD, 2009, pp. 1--15.
\newline\urlprefix\url{http://issfd.org/ISSFD_2009/InterMissionDesignII/Lara.pdf}

\bibitem{Lara2011}
M.~{Lara},
  \href{http://www.sciencedirect.com/science/article/pii/S009457651100066X}{{Design
  of long-lifetime lunar orbits: A hybrid approach}}, Acta Astronautica
  69~(3--4) (2011) 186--199.
\newblock \href {http://dx.doi.org/10.1016/j.actaastro.2011.03.009}
  {\path{doi:10.1016/j.actaastro.2011.03.009}}.
\newline\urlprefix\url{http://www.sciencedirect.com/science/article/pii/S009457651100066X}

\bibitem{Cook1966}
G.~E. {Cook}, {Perturbations of near-circular orbits by the Earth's
  gravitational potential}, Planetary and Space Science 14 (1966) 433--444.
\newblock \href {http://dx.doi.org/10.1016/0032-0633(66)90015-8}
  {\path{doi:10.1016/0032-0633(66)90015-8}}.

\bibitem{RosboroughOcampo1991}
G.~W. Rosborough, C.~Ocampo, {Influence of Higher Degree Zonals on the Frozen
  Orbit Geometry (AAS 91-428)}, in: B.~Kaufman, K.~T. Alfriend, R.~L. Roehrich,
  R.~R. Dasenbrock (Eds.), Astrodynamics 1991, Vol.~76 of Advances in the
  Astronautical Sciences, American Astronautical Society, Univelt, Inc.,
  P.O.~Box 28130, San Diego, California 92198, USA, 1992, pp. 1291--1304.

\bibitem{Cook1992}
R.~A. {Cook}, {The long-term behavior of near-circular orbits in a zonal
  gravity field (AAS 91-463)}, in: B.~Kaufman, K.~T. Alfriend, R.~L. Roehrich,
  R.~R. Dasenbrock (Eds.), Astrodynamics 1991, Vol.~76 of Advances in the
  Astronautical Sciences, American Astronautical Society., Univelt, Inc.,
  P.O.~Box 28130, San Diego, California 92198, USA, 1992, pp. 2205--2221.

\bibitem{Lara2010ICATT}
M.~Lara, {A \textit{Mathematica}\raisebox{.6ex}{\tiny\copyright}--based
  approach to the frozen orbits problem about arbitrary bodies. The case of a
  Lunar orbiter}, in: Astrodynamics Beyond Borders, Proceedings of the 4th
  International Conference on Astrodynamics Tools and Techniques, ESA
  publication WPP-308, ICATT, 2010, pp. 1--8.

\bibitem{CoffeyNealSegermanTravisano1995}
S.~L. Coffey, H.~L. Neal, A.~M. Segerman, J.~J. Travisano, An analytic orbit
  propagation program for satellite catalog maintenance, in: K.~T. Alfriend,
  I.~M. Ross, A.~K. Misra, C.~F. Peters (Eds.), AAS/AIAA Astrodynamics
  Conference 1995, Vol.~90 of Advances in the Astronautical Sciences, American
  Astronautical Society, Univelt, Inc., P.O.~Box 28130, San Diego, California
  92198, USA, 1996, pp. 1869--1892.

\bibitem{Knuth1997}
D.~E. Knuth, {The Art of Computer Programming}, Addison--Wesley--Longman,
  Reading, Massachusetts, 1997.

\bibitem{CoffeyDeprit1980}
S.~{Coffey}, A.~{Deprit}, {Fast evaluation of Fourier series}, {Astronomy and
  Astrophysics} 81 (1980) 310--315.

\bibitem{Deprit1982}
A.~{Deprit}, {Delaunay normalisations}, Celestial Mechanics 26 (1982) 9--21.
\newblock \href {http://dx.doi.org/10.1007/BF01233178}
  {\path{doi:10.1007/BF01233178}}.

\bibitem{HealyTravisano1998}
L.~M. {Healy}, J.~J. {Travisano}, {Automatic rendering of astrodynamics
  expressions for efficient evaluation}, Journal of the Astronautical Sciences
  46~(1) (1998) 65--81.

\bibitem{WnukBreiter1991}
E.~Wnuk, S.~Breiter,
  \href{http://www.sciencedirect.com/science/article/pii/027311779190251E}{{The
  motion of natural and artificial satellites in Mars gravity field}}, Advances
  in Space Research 11~(6) (1991) 183 -- 188.
\newblock \href
  {http://dx.doi.org/https://doi.org/10.1016/0273-1177(91)90251-E}
  {\path{doi:https://doi.org/10.1016/0273-1177(91)90251-E}}.
\newline\urlprefix\url{http://www.sciencedirect.com/science/article/pii/027311779190251E}

\bibitem{Metrisetal1993}
G.~{Metris}, P.~{Exertier}, Y.~{Boudon}, F.~{Barlier}, {Long period variations
  of the motion of a satellite due to non-resonant tesseral harmonics of a
  gravity}, Celestial Mechanics and Dynamical Astronomy 57 (1993) 175--188.
\newblock \href {http://dx.doi.org/10.1007/BF00692472}
  {\path{doi:10.1007/BF00692472}}.

\bibitem{WnukJopek1994}
E.~{Wnuk}, T.~{Jopek}, {Satellite orbit calculations using geopotential
  coefficients up to high degree and order}, Advances in Space Research 14.

\bibitem{LaraPalacianRussell2010}
M.~Lara, J.~Palaci\'an, R.~Russell,
  \href{https://link.springer.com/article/10.1007/s10569-010-9286-2}{{Mission
  design through averaging of perturbed Keplerian systems: the paradigm of an
  Enceladus orbiter}}, Celestial Mechanics and Dynamical Astronomy 108~(1)
  (2010) 1--22.
\newblock \href {http://dx.doi.org/10.1007/s10569-010-9286-2}
  {\path{doi:10.1007/s10569-010-9286-2}}.
\newline\urlprefix\url{https://link.springer.com/article/10.1007/s10569-010-9286-2}

\bibitem{LaraPerezLopez2017}
M.~Lara, I.~P\'erez, R.~L\'opez,
  \href{http://www.sciencedirect.com/science/article/pii/S100757041730429X}{{Higher
  Order Approximation to the Hill Problem Dynamics about the Libration
  Points}}, Communications in Nonlinear Science and Numerical Simulation in
  press.
\newblock \href {http://dx.doi.org/10.1016/j.cnsns.2017.12.007}
  {\path{doi:10.1016/j.cnsns.2017.12.007}}.
\newline\urlprefix\url{http://www.sciencedirect.com/science/article/pii/S100757041730429X}

\bibitem{LaraFerrerSaedeleer2009}
M.~Lara, S.~Ferrer, B.~d. Saedeleer,
  \href{http://dx.doi.org/10.1007/BF03321517}{Lunar analytical theory for polar
  orbits in a 50-degree zonal model plus third-body effect}, The Journal of the
  Astronautical Sciences 57~(3) (2009) 561--577.
\newblock \href {http://dx.doi.org/10.1007/BF03321517}
  {\path{doi:10.1007/BF03321517}}.
\newline\urlprefix\url{http://dx.doi.org/10.1007/BF03321517}

\bibitem{Lara2018Stardust}
M.~{Lara}, {Exploring Sensitivity of Orbital Dynamics with Respect to Model
  Truncation: The Frozen Orbits Approach}, in: M.~Vasile, E.~Minisci,
  L.~Summerer, P.~McGinty (Eds.), Stardust Final Conference, Vol.~52 of
  {Astrophysics and Space Science Proceedings}, Springer International
  Publishing, Cham, 2018, pp. 69--83.
\newblock \href {http://dx.doi.org/10.1007/978-3-319-69956-1_4}
  {\path{doi:10.1007/978-3-319-69956-1_4}}.

\bibitem{Cunningham1970}
L.~E. {Cunningham}, {On the Computation of the Spherical Harmonic Terms Needed
  during the Numerical Integration of the Orbital Motion of an Artificial
  Satellite}, Celestial Mechanics 2 (1970) 207--216.
\newblock \href {http://dx.doi.org/10.1007/BF01229495}
  {\path{doi:10.1007/BF01229495}}.

\bibitem{Pines1973}
S.~{Pines}, {Uniform Representation of the Gravitational Potential and its
  Derivatives}, AIAA Journal 11 (1973) 1508--1511.
\newblock \href {http://dx.doi.org/10.2514/3.50619}
  {\path{doi:10.2514/3.50619}}.

\bibitem{LundbergSchutz1988}
J.~B. {Lundberg}, B.~E. {Schutz}, {Recursion formulas of Legendre functions for
  use with nonsingular geopotential models.}, Journal of Guidance Control
  Dynamics 11 (1988) 31--38.
\newblock \href {http://dx.doi.org/10.2514/3.20266}
  {\path{doi:10.2514/3.20266}}.

\bibitem{FantinoCasotto2009}
E.~{Fantino}, S.~{Casotto}, {Methods of harmonic synthesis for global
  geopotential models and their first-, second- and third-order gradients},
  Journal of Geodesy 83 (2009) 595--619.
\newblock \href {http://dx.doi.org/10.1007/s00190-008-0275-0}
  {\path{doi:10.1007/s00190-008-0275-0}}.

\bibitem{Kozai1959}
Y.~{Kozai}, {The Motion of a Close Earth Satellite}, The Astronomical Journal
  64~(11) (1959) 367--377.

\bibitem{Lara2008}
M.~{Lara}, {Simplified Equations for Computing Science Orbits Around Planetary
  Satellites}, Journal of Guidance Control Dynamics 31~(1) (2008) 172--181.
\newblock \href {http://dx.doi.org/10.2514/1.31107}
  {\path{doi:10.2514/1.31107}}.

\bibitem{LaraPalacianYanguasCorral2010}
M.~{Lara}, J.~F. {Palaci{\'a}n}, P.~{Yanguas}, C.~{Corral},
  \href{http://www.sciencedirect.com/science/article/pii/S0094576509004974}{{Analytical
  theory for spacecraft motion about Mercury}}, Acta Astronautica 66~(7-8)
  (2010) 1022--1038.
\newblock \href {http://dx.doi.org/10.1016/j.actaastro.2009.10.011}
  {\path{doi:10.1016/j.actaastro.2009.10.011}}.
\newline\urlprefix\url{http://www.sciencedirect.com/science/article/pii/S0094576509004974}

\bibitem{Kaula1961}
W.~M. {Kaula}, {Analysis of Gravitational and Geometric Aspects of Geodetic
  Utilization of Satellites}, Geophysical Journal 5 (1961) 104--133.
\newblock \href {http://dx.doi.org/10.1111/j.1365-246X.1961.tb00417.x}
  {\path{doi:10.1111/j.1365-246X.1961.tb00417.x}}.

\bibitem{Kaula1966}
W.~M. {Kaula}, {Theory of satellite geodesy. Applications of satellites to
  geodesy}, Blaisdell, Waltham, Massachusetts, 1966.

\bibitem{Hansen1855}
P.~A. {Hansen}, {Expansions of the product of a power of the radius vector with
  the sinus or cosinus of a multiple of the true anomaly in terms of series
  containing the sinuses or cosinuses of the multiples of the true, eccentric
  or mean anomaly}, Abhandlungen der Koniglich Sachsischen Gesellschaft der
  Wissenschaften 2~(3) (1855) 183--281, {English translation by J.C. Van der
  Ha, ESA/ESOC, Darmstadt, Germany, 1977}.

\bibitem{ArsenaultFordKoskela1970}
J.~L. {Arsenault}, K.~C. {Ford}, P.~E. {Koskela}, {Orbit determination using
  analytic partial derivatives of perturbed motion.}, AIAA Journal 8 (1970)
  4--12.
\newblock \href {http://dx.doi.org/10.2514/3.5597} {\path{doi:10.2514/3.5597}}.

\bibitem{BrouckeCefola1972}
R.~A. {Broucke}, P.~J. {Cefola}, {On the Equinoctial Orbit Elements}, Celestial
  Mechanics 5~(3) (1972) 303--310.
\newblock \href {http://dx.doi.org/10.1007/BF01228432}
  {\path{doi:10.1007/BF01228432}}.

\bibitem{CefolaBroucke1975}
P.~J. {Cefola}, R.~{Broucke}, {On the formulation of the gravitational
  potential in terms of equinoctial variables}, in: 13th Aerospace Sciences
  Meeting, Pasadena, California, American Institute of Aeronautics and
  Astronautics, USA, 1975, pp. 1--25, {AIAA Paper No.~75-9}.

\bibitem{Saedeleer2005}
B.~{de Saedeleer}, {Complete Zonal Problem of the Artificial Satellite: Generic
  Compact Analytic First Order in Closed Form}, Celestial Mechanics and
  Dynamical Astronomy 91 (2005) 239--268.
\newblock \href {http://dx.doi.org/10.1007/s10569-004-1813-6}
  {\path{doi:10.1007/s10569-004-1813-6}}.

\bibitem{Kaula1962}
W.~M. {Kaula}, {Development of the lunar and solar disturbing functions for a
  close satellite}, The Astronomical Journal 67 (1962) 300.
\newblock \href {http://dx.doi.org/10.1086/108729} {\path{doi:10.1086/108729}}.

\bibitem{LaskarBoue2010}
{Laskar, J.}, {Bou\'e, G.},
  \href{https://doi.org/10.1051/0004-6361/201014496}{Explicit expansion of the
  three-body disturbing function for arbitrary eccentricities and
  inclinations}, {Astronomy and Astrophysics} 522 (2010) A60.
\newblock \href {http://dx.doi.org/10.1051/0004-6361/201014496}
  {\path{doi:10.1051/0004-6361/201014496}}.
\newline\urlprefix\url{https://doi.org/10.1051/0004-6361/201014496}

\bibitem{PalacianVanegasYanguas2017}
J.~F. {Palaci{\'a}n}, J.~{Vanegas}, P.~{Yanguas}, {Compact normalisations in
  the elliptic restricted three body problem}, Astrophysics and Space Science
  362 (2017) 215.
\newblock \href {http://dx.doi.org/10.1007/s10509-017-3195-8}
  {\path{doi:10.1007/s10509-017-3195-8}}.

\bibitem{Brouwer1959}
D.~{Brouwer}, {Solution of the problem of artificial satellite theory without
  drag}, The Astronomical Journal 64 (1959) 378--397.
\newblock \href {http://dx.doi.org/10.1086/107958} {\path{doi:10.1086/107958}}.

\bibitem{Garfinkel1959}
B.~{Garfinkel}, {The orbit of a satellite of an oblate planet}, The
  Astronomical Journal 64~(9) (1959) 353--367.
\newblock \href {http://dx.doi.org/10.1086/107956} {\path{doi:10.1086/107956}}.

\bibitem{Deprit1981}
A.~Deprit, The elimination of the parallax in satellite theory, Celestial
  Mechanics 24~(2) (1981) 111--153.
\newblock \href {http://dx.doi.org/10.1007/BF01229192}
  {\path{doi:10.1007/BF01229192}}.

\bibitem{LaraSanJuanLopezOchoa2013b}
M.~{Lara}, J.~F. {San-Juan}, L.~M. {L{\'o}pez-Ochoa}, {Proper Averaging Via
  Parallax Elimination (AAS 13-722)}, in: S.~B. Broschart, J.~D. Turner, K.~C.
  Howell, F.~R. Hoots (Eds.), Astrodynamics 2013, Vol. 150 of {Advances in the
  Astronautical Sciences}, American Astronautical Society, Univelt, Inc.,
  P.O.~Box 28130, San Diego, California 92198, USA, 2014, pp. 315--331.

\bibitem{LaraSanJuanLopezOchoa2013c}
M.~{Lara}, J.~F. {San-Juan}, L.~M. {L{\'o}pez-Ochoa}, {Delaunay variables
  approach to the elimination of the perigee in Artificial Satellite Theory},
  Celestial Mechanics and Dynamical Astronomy 120~(1) (2014) 39--56.
\newblock \href {http://arxiv.org/abs/1312.7577} {\path{arXiv:1312.7577}},
  \href {http://dx.doi.org/10.1007/s10569-014-9559-2}
  {\path{doi:10.1007/s10569-014-9559-2}}.

\bibitem{Deprit1969}
A.~{Deprit}, Canonical transformations depending on a small parameter,
  Celestial Mechanics 1~(1) (1969) 12--30.
\newblock \href {http://dx.doi.org/10.1007/BF01230629}
  {\path{doi:10.1007/BF01230629}}.

\bibitem{BoccalettiPucacco1998v2}
D.~{Boccaletti}, G.~{Pucacco}, {Theory of orbits. Volume 2: Perturbative and
  geometrical methods}, 1st Edition, Astronomy and Astrophysics Library,
  Springer-Verlag, Berlin Heidelberg New York, 2002.

\bibitem{MeyerHall1992}
K.~R. {Meyer}, G.~R. {Hall}, {Introduction to Hamiltonian Dynamical Systems and
  the N-Body Problem}, Springer, New York, 1992.

\bibitem{MacMillan1958}
W.~D. MacMillan, {The Theory of the Potential}, Dover Publishers Inc., New
  York, 1958.

\bibitem{WGS84}
{Defense Mapping Agency}, {Department of Defense World Geodetic System 1984:
  Its definition and relationship with local geodetic systems}, DMA Technical
  Report 8350.2, Centre National d'\'Etudes Spatiales, Washington, D.C. (May
  1987).

\bibitem{KonoplivBanerdtSjogren1999}
A.~Konopliv, W.~Banerdt, W.~Sjogren,
  \href{http://www.sciencedirect.com/science/article/pii/S0019103599960864}{Venus
  gravity: 180th degree and order model}, Icarus 139~(1) (1999) 3 -- 18.
\newblock \href {http://dx.doi.org/https://doi.org/10.1006/icar.1999.6086}
  {\path{doi:https://doi.org/10.1006/icar.1999.6086}}.
\newline\urlprefix\url{http://www.sciencedirect.com/science/article/pii/S0019103599960864}

\bibitem{KonoplivParkFolkner2016}
A.~S. Konopliv, R.~S. Park, W.~M. Folkner,
  \href{http://www.sciencedirect.com/science/article/pii/S0019103516001305}{An
  improved jpl mars gravity field and orientation from mars orbiter and lander
  tracking data}, Icarus 274 (2016) 253 -- 260.
\newblock \href
  {http://dx.doi.org/https://doi.org/10.1016/j.icarus.2016.02.052}
  {\path{doi:https://doi.org/10.1016/j.icarus.2016.02.052}}.
\newline\urlprefix\url{http://www.sciencedirect.com/science/article/pii/S0019103516001305}

\bibitem{Battin1999}
R.~H. Battin, An Introduction to the Mathematics and Methods of Astrodynamics,
  American Institute of Aeronautics and Astronautics, Reston, VA, 1999.

\bibitem{Nayfeh2004}
A.~H. Nay\-feh, Perturbation Methods, Wiley-VCH Verlag GmbH \& Co.~KGaA,
  Weinheim, Germany, 2004.

\bibitem{BrouwerClemence1961}
D.~Brouwer, G.~M. Clemence, Methods of Celestial Mechanics, Academic Press, New
  York and London, 1961.

\bibitem{FerrazMello2007}
S.~{Ferraz-Mello}, {Canonical Perturbation Theories - Degenerate Systems and
  Resonance}, Vol. 345 of Astrophysics and Space Science Library, Springer, New
  York, 2007.

\bibitem{Delaunay1860}
C.~E. {Delaunay}, \href{http://gallica.bnf.fr/ark:/12148/cb343783130/date}{{La
  Th\'eorie du Mouvement de la Lune, Premier volume}}, Vol.~28 of M\'emoires de
  l'Academie des Sciences de l'Institut Imp\'erial de France.,
  {Mallet-Bachellier}, Paris, 1860.
\newline\urlprefix\url{http://gallica.bnf.fr/ark:/12148/cb343783130/date}

\bibitem{Kozai1962}
Y.~{Kozai}, {Second-Order Solution of Artificial Satellite Theory without Air
  Drag}, The Astronomical Journal 67~(7) (1962) 446--461.

\bibitem{OsacarPalacian1994}
C.~{Os\'acar}, J.~F. {Palaci\'an}, {Decomposition of functions for elliptical
  orbits}, Celestial Mechanics and Dynamical Astronomy 60~(2) (1994) 207--223.
\newblock \href {http://dx.doi.org/10.1007/BF00693322}
  {\path{doi:10.1007/BF00693322}}.

\bibitem{DepritRom1970}
A.~Deprit, A.~Rom, {The Main Problem of Artificial Satellite Theory for Small
  and Moderate Eccentricities}, Celestial Mechanics 2~(2) (1970) 166--206.

\bibitem{Wnuk1997}
E.~{Wnuk}, {Highly eccentric satellite orbits}, Advances in Space Research 19
  (1997) 1735--1740.
\newblock \href {http://dx.doi.org/10.1016/S0273-1177(97)00336-0}
  {\path{doi:10.1016/S0273-1177(97)00336-0}}.

\bibitem{BrumbergFukushima1994}
E.~{Brumberg}, T.~{Fukushima}, {Expansions of elliptic motion based on elliptic
  function theory}, Celestial Mechanics and Dynamical Astronomy 60 (1994)
  69--89.
\newblock \href {http://dx.doi.org/10.1007/BF00693093}
  {\path{doi:10.1007/BF00693093}}.

\bibitem{Brumbergetal1995}
E.~{Brumberg}, V.~A. {Brumberg}, T.~{Konrad}, M.~{Soffel}, {Analytical Linear
  Perturbation Theory for Highly Eccentric Satellite Orbits}, Celestial
  Mechanics and Dynamical Astronomy 61 (1995) 369--387.
\newblock \href {http://dx.doi.org/10.1007/BF00049516}
  {\path{doi:10.1007/BF00049516}}.

\bibitem{Fukushima2014ch}
T.~Fukushima, Elliptic functions and elliptic integrals for celestial mechanics
  and dynamical astronomy, in: S.~Kopeikin (Ed.), Frontiers in Relativistic
  Celestial Mechanics Volume 2: Applications and Experiments, Vol.~22 of De
  Gruyter Studies in Mathematical Physics, Walter de Gruyter GmbH, Genthiner
  Strasse 13, D-10785 Berlin / Germany, 2014, pp. 189--228.

\bibitem{AlfriendCoffey1984}
K.~T. {Alfriend}, S.~L. {Coffey}, {Elimination of the perigee in the satellite
  problem}, Celestial Mechanics 32~(2) (1984) 163--172.
\newblock \href {http://dx.doi.org/10.1007/BF01231123}
  {\path{doi:10.1007/BF01231123}}.

\bibitem{ExertierThesis1988}
A.~Exertier, Orbitographie des satellites artificiels sur de grandes periodes
  de temps. Possibilites d'applications, PhD.~Thesis, Observatoire de Paris,
  Paris, 1988.

\bibitem{LaraSanJuanLopezCefola2012}
M.~{Lara}, J.~F. {San-Juan}, L.~M. {L{\'o}pez}, P.~J. {Cefola}, {On the
  third-body perturbations of high-altitude orbits}, Celestial Mechanics and
  Dynamical Astronomy 113 (2012) 435--452.
\newblock \href {http://dx.doi.org/10.1007/s10569-012-9433-z}
  {\path{doi:10.1007/s10569-012-9433-z}}.

\bibitem{CoffeyAlfriend1981}
S.~Coffey, K.~T. Alfriend, Short period elimination for the tesseral harmonics,
  in: A.~L. Friedlander, P.~J. Cefola, B.~Kaufman, W.~Williamson, G.~Tseng
  (Eds.), AAS/AIAA Astrodynamics Conference 1981, Vol.~46 of Advances in the
  Astronautical Sciences, American Astronautical Society, Univelt, Inc.,
  P.O.~Box 28130, San Diego, California 92198, USA, 1982, pp. 87--101.

\bibitem{CoffeyDeprit1982}
S.~Coffey, A.~Deprit, {Third-Order Solution to the Main Problem in Satellite
  Theory}, Journal of Guidance, Control and Dynamics 5~(4) (1982) 366--371.

\bibitem{DanbyDepritRom1965}
J.~M.~A. {Danby}, A.~{Deprit}, A.~R.~M. {Rom}, {The Symbolic Manipulation of
  Poisson Series}, {Mathematical Note No. 432} D1-82-0481, Boeing Scientific
  Research Laboratories, Seattle, Washington (1965).

\bibitem{DepritRom1968}
A.~{Deprit}, A.~{Rom}, {Lindstedt's Series on a Computer}, The Astronomical
  Journal 73 (1968) 210.
\newblock \href {http://dx.doi.org/10.1086/110619} {\path{doi:10.1086/110619}}.

\bibitem{Lara2010}
M.~{Lara},
  \href{http://www.aeroespacial.org.br/jaesa/editions/repository/v02/n01/6-Lara.pdf}{{Three-body
  dynamics around the smaller primary. Application to the design of science
  orbits}}, Journal of Aerospace Engineering, Sciences and Applications 2~(1)
  (2010) 53--65.
\newblock \href {http://dx.doi.org/10.7446/jaesa.0201.06}
  {\path{doi:10.7446/jaesa.0201.06}}.
\newline\urlprefix\url{http://www.aeroespacial.org.br/jaesa/editions/repository/v02/n01/6-Lara.pdf}

\bibitem{LaraSanJuanLopezOchoa2013a}
M.~{Lara}, J.~F. {San-Juan}, L.~M. {L{\'o}pez-Ochoa}, {Precise analytical
  computation of frozen-eccentricity, low Earth orbits in a tesseral
  potential}, Mathematical Problems in Engineering 2013~(Article ID 191384)
  (2013) 1--13.
\newblock \href {http://dx.doi.org/10.1155/2013/ 191384}
  {\path{doi:10.1155/2013/ 191384}}.

\bibitem{Shapiro1996}
B.~E. Shapiro, {Phase Plane Analysis and Observed Frozen Orbit for the
  Topex/Poseidon Mission}, in: P.~M. Bainum, G.~L. May, Y.~Ohkami, K.~Uesugi,
  Q.~Faren, L.~Furong (Eds.), 6th AAS/JRS/CSA Symposium, International Space
  Conference of Pacific-Basin Societies, Strengthening Cooperation in the 21st
  Century, Vol.~91 of Advances in the Astronautical Sciences, American
  Astronautical Society, Univelt, Inc., P.O.~Box 28130, San Diego, California
  92198, USA, 1996, pp. 853--872.

\bibitem{FoltaQuinn2006}
D.~{Folta}, D.~{Quinn}, {Lunar Frozen Orbits}, in: AIAA/AAS Astrodynamics
  Specialist Conference and Exhibit, 21 - 24 August 2006, Keystone, Colorado,
  American Institute of Aeronautics and Astronautics, 2006, pp. 1--18, {AIAA
  Paper 2006-6749}.
\newblock \href {http://dx.doi.org/10.2514/6.2006-6749}
  {\path{doi:10.2514/6.2006-6749}}.

\bibitem{Konoplivetal2001}
A.~S. {Konopliv}, S.~W. {Asmar}, E.~{Carranza}, W.~L. {Sjogren}, D.~N. {Yuan},
  {Recent Gravity Models as a Result of the Lunar Prospector Mission}, Icarus
  150 (2001) 1--18.
\newblock \href {http://dx.doi.org/10.1006/icar.2000.6573}
  {\path{doi:10.1006/icar.2000.6573}}.

\end{thebibliography}
\end{document}